\documentclass[10pt,a4paper]{article}
\usepackage[utf8]{inputenc}
\usepackage[english]{babel}
\usepackage{amsmath}
\usepackage{amsfonts}
\usepackage{amssymb}
\usepackage{mathtools}
\usepackage{amsthm}
\usepackage{float}
\usepackage{etoolbox}
\usepackage{bbm}
\usepackage{hyperref}
\usepackage{graphicx}
\usepackage[table]{xcolor}
\usepackage[T1]{fontenc}
\newtheorem{theorem}{Theorem}
\newtheorem{lemma}{Lemma}
\newtheorem{remark}{Remark}
\newtheorem{definition}{Definition}
\newtheorem{corollary}{Corollary}

\usepackage[left=2cm,right=2cm,top=2cm,bottom=2cm]{geometry}

\newcommand{\probability}{\operatorname{\mathbb{P}}\probarg}
\DeclarePairedDelimiterX{\probarg}[1]{(}{)}{%
  \ifnum\currentgrouptype=16 \else\begingroup\fi
  \activatebar#1
  \ifnum\currentgrouptype=16 \else\endgroup\fi
}
\newcommand{\activatebar}{%
  \begingroup\lccode`\~=`\|
  \lowercase{\endgroup\let~}\innermid 
  \mathcode`|=\string"8000
}
\newcommand{\LimitN}{\overset{N\to\infty}{\longrightarrow}}
\newcommand{\LimitD}{\overset{d}{\longrightarrow}}

\newcommand{\TailExpon}[1]{\gamma(#1)}
\newcommand{\TailEvent}[1]{f_N(#1)}
\newcommand{\Process}[2]{W_{#1}(#2)+W_A(#2)-\beta #2}

\newcommand{\HTime}[3]{\tau_{#1,N}^{#2,#3}}
\newcommand{\CondProb}[2]{\operatorname{\mathbb{P}}^{(#1<#2)}\probarg}
\newcommand{\HTimeA}[2]{\tilde{\tau}^{#1,#2}_{A,N}}

\newcommand{\BrQueueLength}[3]{\ifstrequal{#2}{}{Q^{\beta}_{#1,A}}{\ifstrequal{#3}{}{Q^{\beta}_{#1,A}(#2,\infty)}{Q^{\beta}_{#1,A}(#2,#3)}}}

\newcommand{\HatBrQueueLength}[3]{\ifstrequal{#2}{}{\hat{Q}^{\beta}_{#1,A}}{\ifstrequal{#3}{}{\hat{Q}^{\beta}_{#1,A}(#2,\infty)}{\hat{Q}^{\beta}_{#1,A}(#2,#3)}}}

\newcommand{\MaxQueue}[2]{\ifstrequal{#1}{}{\bar Q_N^{\beta}}{\ifstrequal{#2}{}{\bar Q_N^{\beta}(#1,\infty)}{\bar Q_N^{\beta}(#1,#2)}}}

\author{

Dennis Schol\\
\small{Eindhoven University of Technology, P.O. Box 513, 5600 MB Eindhoven,}\\ \small{The Netherlands, c.schol@tue.nl}\\ \\
Maria Vlasiou\\ 
\small{University of Twente, Eindhoven University of Technology, P.O. Box 513, 5600 MB Eindhoven,}\\ \small{The Netherlands, m.vlasiou@tue.nl}\\ \\
Bert Zwart\\ 
\small{Eindhoven University of Technology, CWI, P.O. Box 513, 5600 MB Eindhoven,}\\ \small{The Netherlands,  bert.zwart@cwi.nl}\\
}

\hypersetup{%
pdfpagemode=UseNone,
pdftitle={},%
pdfauthor={Dennis},%
pdfkeywords={paper3},%
pdfstartview={FitH},%
colorlinks=true,%
citecolor=black,%
linkcolor=black,%
urlcolor=black,%
unicode=true,%
breaklinks=true%
}%

\begin{document}
\title{Tail Asymptotics for the Delay in a Brownian Fork-Join Queue}
\maketitle
\begin{abstract}
    In this paper, we study the tail behavior of  $\max_{i\leq N}\sup_{s>0}\left(\Process{i}{s}\right)$ as $N\to\infty$, with $(W_i,i\leq N)$ i.i.d.\ Brownian motions and $W_A$ an independent Brownian motion. This random variable can be seen as the maximum of $N$ mutually dependent Brownian queues, which in turn can be interpreted as the backlog in a Brownian fork-join queue. In previous work, we have shown that this random variable centers around $\frac{\sigma^2}{2\beta}\log N$. Here, we analyze the rare-event that this random variable reaches the value $(\frac{\sigma^2}{2\beta}+a)\log N$, with $a>0$. It turns out that its probability behaves roughly as a power law with $N$, where the exponent depends on $a$. However, there are three regimes, around a critical point $a^{\star}$; namely, $0<a<a^{\star}$, $a=a^{\star}$, and $a>a^{\star}$. The latter regime exhibits a form of asymptotic independence, while the first regime reveals highly irregular behavior with a clear dependence structure among the $N$ suprema, with a nontrivial transition at $a=a^{\star}$.
\end{abstract}
Keywords: Brownian queues; fork-join queues; extreme value theory; tail asymptotics

\section{Introduction}
Fork-join queues are a useful modeling tool for congestion in complex networks, such as assembly systems, communication networks, and supply chains. Such networks can be large and assembly is only possible upon availability of all parts. Thus, the bottleneck of the system is caused by the slowest production line in the system. This setting motivates us to investigate such delays in a stylized version of a large fork-join queueing system. In this setting, a key quantity of interest is the behavior of the longest queue. We assume that arrival and service processes are Brownian, as it is a standard result in queueing theory that queueing systems in heavy-traffic can be approximated by reflected Brownian motions. Furthermore, when the arrival and service processes are deterministic with some white noise perturbation, it is also a natural choice to model this with Brownian motions. We analyze the steady-state behavior of this system. Hence, we can model the backlog in queue $i$ by
$\BrQueueLength{i}{}{}=\sup_{s>0}(W_i(s)+W_A(s)-\beta s)$, where $W_A$ is a Brownian motion term with standard deviation $\sigma_A$ that represents the fluctuations in the arrival process, $W_i$ is a Brownian motion term with standard deviation $\sigma$ that represents the fluctuations in the service process, and $\beta>0$ represents the drift of the queue. Furthermore, we assume that $(W_i,i\leq N)$ are i.i.d.\ Brownian motions, and for all $i$, $W_i$ and $W_A$ are mutually independent. These are natural choices as well, because these assumptions indicate that servers' work speeds are mutually independent, and independent with respect to the interarrival times.

Because the bottleneck in the system is the slowest production line, we are interested in the longest queue length, and we investigate the random variable
$\MaxQueue{}{}=\max_{i\leq N}\BrQueueLength{i}{}{}$. We see that this random variable is a maximum of $N$ dependent random variables, due to the common arrival process $W_A$. As we try to model systems with many servers, we are typically interested in the behavior of this random variable as $N\to\infty$. In \cite{MeijerSchol2021}, it is shown that $\MaxQueue{}{}$ is in the domain of attraction of the normal distribution:
\begin{align}\label{eq: central limit convergence}
    \probability*{\MaxQueue{}{}>\frac{\sigma^2}{2\beta}\log N+x\sqrt{\log N}}\LimitN\probability*{\frac{\sigma\sigma_A}{\sqrt{2}\beta}X>x},
\end{align}
with $X\overset{d}{=}\mathcal{N}(0,1)$. This means that $\MaxQueue{}{}$ centers around $\frac{\sigma^2}{2\beta}\log N$ and deviates with order $\sqrt{\log N}$.

This convergence result provides a prediction of the typical delay. However, one might also be interested in the question how likely it is that the delay will be much longer, as delays may cause large costs. Obviously, the probability $\probability{\MaxQueue{}{}>y_N}\LimitN 0$, when $y_N-\frac{\sigma^2}{2\beta}\log N$ grows to infinity at a rate faster than $\sqrt{ \log N}$, but the question is how fast this probability converges to 0. In this study, we focus on the probability
\begin{align*}
     \probability*{\MaxQueue{}{}>\left(\frac{\sigma^2}{2\beta}+a\right)\log N},
\end{align*}
with $a>0$. As we show later on, the exact behavior of this tail probability depends on the choice of $a$, where we can distinguish three regimes: $0<a<a^{\star}$, $a=a^{\star}$, and $a>a^{\star}$, with $a^{\star}$ an explicitly identified constant in $(0,\infty)$. The logarithmic asymptotics for these three regimes are given in Theorem \ref{thm: log asymp}, while the sharper asymptotics for the cases $a>a^{\star}$, $a=a^{\star}$, and $0<a<a^{\star}$ are given in Theorems \ref{thm: exact asymp a > a star}, \ref{thm: exact asymp a = a star}, and \ref{thm: exact asymp 0<a< a star}, respectively. It easily follows from the proofs that when $y_N$ is of larger order than $\log N$, the convergence behavior of $\probability{\MaxQueue{}{}>y_N}$ is the same as for the case $a>a^{\star}$, cf.\ Corollary \ref{cor: yN larger than log N}.

Our work is related to the literature on extreme values of Gaussian processes. In this paper, we examine exceedance probabilities of the order $(\frac{\sigma^2}{2\beta}+a)\log N$ with $a>0$. More work has been done on joint suprema of Brownian motions. For instance, \cite{kou2016first} gives the solution of the Laplace transform of joint first passage times in terms of the solution of a partial differential equation, where the Brownian motions are dependent. Further, \cite{debicki2020exact} analyze the tail asymptotics of the all-time suprema of two dependent Brownian motions. The joint suprema of a finite number of Brownian motions is also studied  \cite{debicki2015extremes}, where the authors give tail asymptotics of the joint suprema of independent Gaussian processes over a finite time interval. These are just three examples -- more results may be found in  \cite{mandjes2007large} and \cite{piterbarg1996asymptotic}. 

Our work also relates to the literature on fork-join queues. Exact results on fork-join queues with two service stations can be found in \cite{baccelli1985two,flatto1984two,de1988fredholm,wright1992two}. Approximations for systems with an arbitrary but fixed number of servers can be found in \cite{baccelli1989queueing,ko2004response,nelson1988approximate}. In \cite{varma1990heavy} a heavy-traffic analysis for fork-join queues is derived; see also \cite{nguyen1993processing} and \cite{nguyen1994trouble}. More recent work in this direction may be found in 
\cite{lu2015gaussian,lu2017heavy,lu2017heavy2,scholMOR}. Our work adds to the existing literature, as we analyze the largest of $N$ queues as $N\to\infty$. Literature on such extreme value results is rare. More specifically, we derive a large deviation principle for the longest of $N$ dependent Brownian queues as $N\to\infty$, to obtain this, we use and extend the results obtained in \cite{debicki2020exact}, in which the case $N=2$ is investigated.

This paper is organized as follows. In Section \ref{sec: main results}, we present our main results, which contain an interesting phase transition in the way a large supremum occurs depending on the value of $a$. We explain the reason behind this phase transition in detail. 
The rest of the paper is devoted to proofs. In Section \ref{sec: log asymp}, we give a proof of Theorem \ref{thm: log asymp}, which focuses on logarithmic asymptotics. In Section \ref{sec: useful lemmas}, we present some auxiliary lemmas. In Sections \ref{subsec: a > a star}, \ref{subsec: a = a star}, and \ref{subsec: 0<a< a star}, we provide the proofs of Theorems \ref{thm: exact asymp a > a star}, \ref{thm: exact asymp a = a star}, and \ref{thm: exact asymp 0<a< a star}, respectively, which deal with asymptotic estimates that are sharper than Theorem \ref{thm: log asymp}.


\section{Main results}\label{sec: main results}
\begin{figure}[H]
\centering
\includegraphics[scale=0.75]{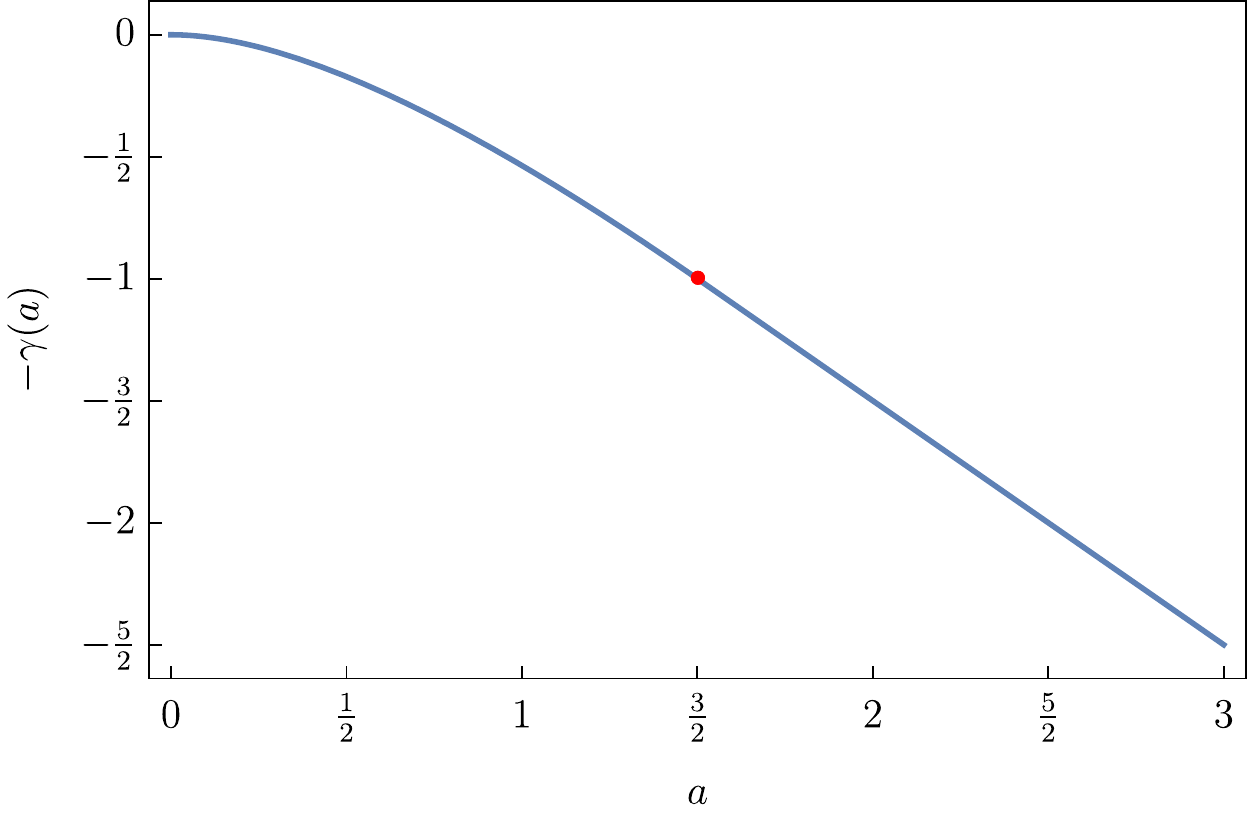}
\caption{$\sigma=1, \sigma_A=1, \beta=1, a^{\star}=3/2$}
\label{fig: tailexpon}
\end{figure}

In this section, we present our main results and also provide some intuition.
We first introduce some additional notation. Recall that $(W_i,i\leq N)$ is a sequence of i.i.d.\ Brownian motions with standard deviation $\sigma$, $W_A$ is a Brownian motion with standard deviation $\sigma_A$, $W_i$ and $W_A$ are mutually independent for all $i$, the steady-state queue length in front of server $i$ is given by
\begin{align}
    \BrQueueLength{i}{}{}=\sup_{s>0}(\Process{i}{s}),
\end{align}
and the maximum queue length equals 
\begin{align}\label{eq: def max queue length}
    \MaxQueue{}{}=\max_{i\leq N}\BrQueueLength{i}{}{}.
\end{align}
Further, we write the supremum of a Brownian motion $\{W_i(t)+W_A(t)-\beta t,t>0\}$ over an interval $(u,v)$ as
\begin{align}\label{eq: queue length}
    \BrQueueLength{i}{u}{v}=\sup_{u<s<v}(\Process{i}{s}),
\end{align}
and the maximum of $N$ of these identically distributed random variables as
\begin{align}\label{eq: queue length finite time}
\MaxQueue{u}{v}=\max_{i\leq N}\BrQueueLength{i}{u}{v}.
\end{align}
Also, we introduce shorthand notation that we use later on:
\begin{align}\label{eq: tailevent}
    \TailEvent{a}=\left(\frac{\sigma^2}{2\beta}+a\right)\log N,\\
    \lambda(a)=1-\sigma/\sqrt{2 a\beta+\sigma^2},\label{subeq: def lambda}\\
    T_N(a,k)=\TailEvent{a}/\beta+k\sqrt{\log N},\label{subeq: TN(a,k) def}\\
    T_N(a)=T_N(a,0).\label{subeq: TN(a) def}
\end{align}
Finally, we write 
\begin{align}\label{eq: tail expon}
      \TailExpon{a}=\begin{cases}
     \displaystyle \frac{2a\beta+2\sigma^2-2\sigma\sqrt{2a\beta+\sigma^2}}{\sigma_A^2}&\quad \text{ if }0<a<a^{\star},\\
      \displaystyle \frac{2a\beta-\sigma_A^2}{\sigma^2+\sigma_A^2}& \quad\text{ if } a\geq a^{\star},
      \end{cases}
\end{align}
with $a^{\star}=\frac{\sigma_A^4}{\sigma^22\beta}+\frac{\sigma_A^2}{\beta}$. The function $\TailExpon{a}$ appears in the limit of the logarithmic asymptotics of $\probability{\MaxQueue{}{}>\TailEvent{a}}$.  In Figure \ref{fig: tailexpon}, we plot $-\TailExpon{a}$ for certain choices of the parameters $\sigma,\sigma_A,\beta$, and $a^{\star}$. As can be seen, from $a=a^{\star}$ onwards, the function is linear. Moreover, we see that $\TailExpon{a}$ is continuous everywhere, also for $a=a^{\star}$. 

Our first result, Theorem \ref{thm: log asymp}, provides the logarithmic asymptotics of the tail probability of the maximum steady-state queue length $\probability{\MaxQueue{}{}>\TailEvent{a}}$.
\begin{theorem}\label{thm: log asymp}
Let $a>0$, $(W_i,i\leq N)$ be i.i.d.\ Brownian motions with standard deviation $\sigma$, $W_A$ be a Brownian motion with standard deviation $\sigma_A$, for all $i$, $W_i$ and $W_A$ are mutually independent, and $\MaxQueue{}{}$, $\TailExpon{a}$, and $\TailEvent{a}$ are given by Equations \eqref{eq: def max queue length}, \eqref{eq: tail expon}, and \eqref{eq: tailevent}, respectively, then
\begin{align}\label{eq: log asymp}
\frac{\log(\probability{\MaxQueue{}{}>\TailEvent{a}})}{\log N}\LimitN -\TailExpon{a}.
\end{align}
\end{theorem}
We give the proof of Theorem \ref{thm: log asymp} in Section \ref{sec: log asymp}. To provide some intuition, 
the form of the function $\TailExpon{a}$ suggests there are at least two regimes: the case where $0<a<a^{\star}$, and the case where $a\geq a^{\star}$. These two cases reveal interesting information on the tail behavior of the maximum queue length $\MaxQueue{}{}$.

\paragraph{Case $a>a^{\star}$.} First of all, observe that for $a>a^{\star}$ and $N$ large, by using the convergence result in \eqref{eq: log asymp} and the memoryless property of the exponential distribution, we have that
\begin{equation}\label{eq: heuristic ineq memoryless}
\begin{aligned}
  &\probability*{\max_{i\leq N}\sup_{s>0}\left(\Process{i}{s}\right)>\TailEvent{a}}\\
   &\quad\geq\probability*{\max_{i\leq N}\sup_{s>0}\left(\Process{i}{s}\right)>\TailEvent{a^{\star}}}\probability*{\sup_{s>0}\left(\Process{i}{s}\right)>(a-a^{\star})\log N}\\
   &\quad= \probability*{\max_{i\leq N}\sup_{s>0}\left(\Process{i}{s}\right)>\TailEvent{a^{\star}}}\exp\left(-\frac{-2\beta(a-a^{\star})}{\sigma^2+\sigma_A^2}\log N\right)\\
   &\quad\approx N^{-\TailExpon{a^{\star}}}\exp\left(-\frac{-2\beta(a-a^{\star})}{\sigma^2+\sigma_A^2}\log N\right)\\
   &\quad=N^{-\TailExpon{a}}.
   \end{aligned}
\end{equation}
To understand the lower bound in this expression, observe that due to the memoryless property of an exponentially distributed random variable $E$, we have  that $\probability{E>x+y}=\probability{E>x}\probability{E>x+y\mid E>x}=\probability{E>x}\probability{E>y}$. Then for a sequence of exponentially and identically distributed random variables $(E_i,i\leq N)$, we have for all $j\leq N$ that $\probability{\max_{i\leq N}E_i>x+y}\geq\probability{\max_{i\leq N}E_i>x}\probability{E_j>x+y\mid E_j>x}=\probability{\max_{i\leq N}E_i>x}\probability{E_j>y}$. So, the fact that the tail probability of the maximum steady-state queue length in \eqref{eq: heuristic ineq memoryless} is bounded from below by the expression in \eqref{eq: heuristic ineq memoryless} implies that for $a>a^{\star}$
$$
\probability*{\#\{j\leq N: \sup_{s>0}(W_j(s)+W_A(s)-\beta s)>\TailEvent{a} \}=1\bigg| \MaxQueue{}{}{}>\TailEvent{a}}\LimitN 1.
$$ 
Second, we see that for $a\geq a^{\star}$, $N^{-\TailExpon{a}}=N\probability{\BrQueueLength{i}{}{}>\TailEvent{a}}$. Obviously, since $a\geq 0$, the union bound gives that
\begin{align}\label{eq: union bound}
\probability{\MaxQueue{}{}>\TailEvent{a}}\leq N\probability{\BrQueueLength{i}{}{}>\TailEvent{a}}=N^{-\frac{2a\beta-\sigma_A^2}{\sigma^2+\sigma_A^2}}.
\end{align}
The fact that the union bound is sharp when $a\geq a^{\star}$ indicates that for $a\geq a^{\star}$ the $N$ queues are asymptotically independent; i.e.,
\begin{align*}
   \probability*{\max_{i\leq N}\sup_{s>0}\left(\Process{i}{s}\right)>\TailEvent{a}} \approx \probability*{\max_{i\leq N}\sup_{s>0}\left(W_i(s)+W_{A,i}(s)-\beta s\right)>\TailEvent{a}},
\end{align*}
where the arrival processes $(W_{A,i},i\leq N)$ are independent Brownian motions. In Section \ref{subsec: a = a star}, we see that the boundary case $a=a^{\star}$ does show some dependent behavior, but this dependence structure cannot be deduced from the logarithmic asymptotics.

\paragraph{Case $0<a<a^{\star}$.} 
Finally, the case $0<a<a^{\star}$ is more involved. The function $\gamma(a)$ involves $a$ in a nonlinear fashion. As we observe in Equation \eqref{eq: union bound}, due to the fact that the exponent of the tail probability of an exponentially distributed random variable is linear in $a$, we expect that the logarithmic asymptotics would also be linear in $a$. Thus the structure of $\TailExpon{a}$ shows that the dependent part $W_A$ influences the tail asymptotics, and contrary to the case where $a>a^{\star}$, we have that
$$
\liminf_{N\to\infty}\probability*{\#\{j\leq N: \sup_{s>0}(W_j(s)+W_A(s)-\beta s)>\TailEvent{a} \}>1\bigg| \MaxQueue{}{}{}>\TailEvent{a}}>0.
$$ 
The reason that we see this is that in order to get that the maximum steady-state queue length $\MaxQueue{}{}{}$ reaches the level $\TailEvent{a}$, the arrival process $\{W_A(t)-\lambda(a)\beta t,t>0\}$ must reach a high level around $\lambda(a)\TailEvent{a}$, which is a rare event. Furthermore, one of the $N$ service processes needs to reach a level around $(1-\lambda(a))\TailEvent{a}$; however, this is not a rare event. Even more, the event that a finite number of service processes reaches a level around $(1-\lambda(a))\TailEvent{a}$ has a finite probability. 

The function $\TailExpon{a}$ has more characteristics that can be explained from \cite{MeijerSchol2021}. What we namely see is that $\TailExpon{0}=0$, which is to be expected as we know from \eqref{eq: central limit convergence} and \eqref{eq: tailevent} that for $x=0$
$$
    \probability{\MaxQueue{}{}>\TailEvent{0}}\LimitN \frac{1}{2}.
$$
We further have that 
$\log N \TailExpon{x/\sqrt{\log N}}\LimitN \frac{x^2\beta^2}{\sigma^2\sigma_A^2}.$
It thus follows that for $N$ large,
$$
    N^{-\TailExpon{x/\sqrt{\log N}}}\approx N^{-\frac{x^2\beta^2}{\sigma^2\sigma_A^2\log N}}=\exp\left(-\frac{x^2\beta^2}{\sigma^2\sigma_A^2}\right),
$$
which is the exponent of the limiting distribution given in \eqref{eq: central limit convergence}.

To prove the logarithmic asymptotics in Theorem \ref{thm: log asymp}, it suffices to look at random variables of the type $\max_{i\leq N}(W_i(T_N)+W_A(T_N)-\beta T_N)$ instead of the random variable $\MaxQueue{}{}{}=\max_{i\leq N}\sup_{s>0}(\Process{i}{s})$, where the appropriate choice of $T_N$ is $T_N(a)$, cf.\ Equation \eqref{subeq: TN(a) def}. We show this in more detail in the proof of Lemma \ref{lem: liminf log asymp}. 
For $a>a^{\star}$, the logarithmic asymptotics are relatively straightforward to derive because we see a notion of asymptotic independence, as explained above. In the proof of Lemma \ref{lem: liminf log asymp}, we show that when $0<a\leq a^{\star}$,
\begin{align}\label{eq: approximate log asymptotics}
\begin{split}
    &\log(\probability{\MaxQueue{}{}>\TailEvent{a}})\\
    &\quad\approx\log(\probability{\max_{i\leq N}W_i\big(T_N(a)\big)-(1-\lambda(a))\beta T_N(a)>(1-\lambda(a))\TailEvent{a}})\\
    &\quad\quad+\log(\probability{W_A\big(T_N(a)\big)-\lambda(a)\beta T_N(a)>\lambda(a)\TailEvent{a}}),
\end{split}
    \end{align}
when $N$ is large, and we show that the term $\log(\probability{\max_{i\leq N}W_i\big(T_N(a)\big)-(1-\lambda(a))\beta T_N(a)>(1-\lambda(a))\TailEvent{a}})$ becomes negligible as $N\rightarrow\infty$.


We now turn to precise asymptotics, which are stated in Theorems \ref{thm: exact asymp a > a star}, \ref{thm: exact asymp a = a star}, and \ref{thm: exact asymp 0<a< a star} below for the cases $a>a^{\star}$, $a=a^{\star}$, and $0<a<a^{\star}$, respectively. 
The proofs of these theorems can be found in  Sections \ref{subsec: a > a star}, \ref{subsec: a = a star}, and \ref{subsec: 0<a< a star}.
\begin{theorem}\label{thm: exact asymp a > a star}
Let $a>a^{\star}$, $(W_i,i\leq N)$ be i.i.d.\ Brownian motions with standard deviation $\sigma$, $W_A$ be a Brownian motion with standard deviation $\sigma_A$, for all $i$, $W_i$ and $W_A$ are mutually independent, and $\MaxQueue{}{}$, $\TailExpon{a}$, and $\TailEvent{a}$ are given by Equations \eqref{eq: def max queue length}, \eqref{eq: tail expon}, and \eqref{eq: tailevent}, respectively, then
\begin{align}\label{eq: a>a star thm}
N^{\TailExpon{a}}\probability{\MaxQueue{}{}>\TailEvent{a}}\LimitN 1.    
\end{align}
\end{theorem}
The theorem shows that for $a>a^{\star}$, the tail probability of the steady-state maximum queue length has the same asymptotic behavior as the one for independently and identically distributed arrival processes for each queue. 
\begin{theorem}\label{thm: exact asymp a = a star}
Let $a=a^{\star}$, $(W_i,i\leq N)$ be i.i.d.\ Brownian motions with standard deviation $\sigma$, $W_A$ be a Brownian motion with standard deviation $\sigma_A$, for all $i$, $W_i$ and $W_A$ are mutually independent, and $\MaxQueue{}{}$, $\TailExpon{a}$, and $\TailEvent{a}$ are given by Equations \eqref{eq: def max queue length}, \eqref{eq: tail expon}, and \eqref{eq: tailevent}, respectively, then
\begin{align}\label{eq: a=a star thm}
N^{\TailExpon{a^{\star}}}\probability{\MaxQueue{}{}>\TailEvent{a^{\star}}}\LimitN \frac{1}{2}.    
\end{align}
\end{theorem}
To give a heuristic explanation why we have a transition point at $a=a^{\star}$, recall that $\lambda(a)$ is given in Equation \eqref{subeq: def lambda}, $W_i$ is a Brownian motion with standard deviation $\sigma$, and $W_A$ is a Brownian motion with standard deviation $\sigma_A$. Because the all-time supremum of a Brownian motion is exponentially distributed it is easy to see that for $a=a^{\star}$,
\begin{align*}
    \sup_{s>0}(W_A(s)-\lambda(a^{\star})\beta s)\overset{d}{=}\sup_{s>0}(W_i(s)-(1-\lambda(a^{\star}))\beta s)\overset{d}{=}\sup_{s>0}(W_i(s)+W_A(s)-\beta s).
\end{align*}
Similarly, after a straightforward calculation we observe that for $0<a<a^{\star}$,
 \begin{align*}
     \sup_{s>0}(W_A(s)-\lambda(a)\beta s)\geq_{st.}\sup_{s>0}(W_i(s)-(1-\lambda(a))\beta s),
 \end{align*}
 and for $a>a^{\star}$,
 \begin{align*}
     \sup_{s>0}(W_A(s)-\lambda(a)\beta s)\leq_{st.}\sup_{s>0}(W_i(s)-(1-\lambda(a))\beta s).
 \end{align*} 
 For $0<a<a^{\star}$, large values of $\MaxQueue{}{}$ are predominantly caused by fluctuations of $\{W_A(t)-\lambda(a)\beta t,t>0\}$; we show this rigorously in Section \ref{subsec: 0<a< a star}. In contrast, for $a>a^{\star}$, fluctuations are caused by a combination of the arrival process and one of the service processes, and therefore we see a notion of asymptotic independence.

To explain in more detail why we have a constant 1/2 at the boundary case $a=a^\star$, we first observe that, since the all-time supremum of a Brownian motion with negative drift is exponentially distributed, $\probability{\sup_{s>0}(W_A(s)-\lambda(a^{\star})\beta s)>\lambda(a^{\star})\TailEvent{a^{\star}}}=N^{-\TailExpon{a^{\star}}}$.
Moreover, if the event $\sup_{s>0}(W_A(s)-\lambda(a^{\star})\beta s)>\lambda(a^{\star})\TailEvent{a^{\star}}$ happens, it most likely occurs at time
$T_N(a^{\star})$. By using the union bound 
and that all suprema are the same in distribution we may therefore write
\begin{align*}
&\probability{\MaxQueue{T_N(a^{\star})}{}>\TailEvent{a^{\star}}\mid W_A\big(T_N(a^{\star})\big)-\lambda(a^{\star})\beta T_N(a^{\star})=\lambda(a^{\star})\TailEvent{a^{\star}}}\nonumber\\
    &\quad=\probability*{\max_{i\leq N}\left(W_i\big(T_N(a^{\star})\big)-(1-\lambda(a^{\star}))\beta T_N(a^{\star})+\HatBrQueueLength{i}{}{}\right)>(1-\lambda(a^{\star}))\TailEvent{a^{\star}}}\nonumber\\
    &\quad\approx N\probability*{W_i\big(T_N(a^{\star})\big)-(1-\lambda(a^{\star}))\beta T_N(a^{\star})+\HatBrQueueLength{i}{}{}>(1-\lambda(a^{\star}))\TailEvent{a^{\star}}}\nonumber\\
    &\quad =N\probability*{\sup_{s>T_N(a^{\star})}(W_i(s)-(1-\lambda(a^{\star}))\beta s)>(1-\lambda(a^{\star}))\TailEvent{a^{\star}}}\LimitN \frac{1}{2}.
\end{align*}
If we condition on $\max_{i\leq N}\sup_{s>0}(W_i(s)-(1-\lambda(a^{\star}))\beta s)=(1-\lambda(a^{\star}))\TailEvent{a^{\star}}$, we obtain the same expression after using the same heuristic argument.

Our final result is an  improvement of the logarithmic asymptotics for the case $0<a<a^{\star}$.
\begin{theorem}\label{thm: exact asymp 0<a< a star}
Let $0<a<a^{\star}$, $(W_i,i\leq N)$ be i.i.d.\ Brownian motions with standard deviation $\sigma$, $W_A$ be a Brownian motion with standard deviation $\sigma_A$, for all $i$, $W_i$ and $W_A$ are mutually independent, and $\MaxQueue{}{}$, $\TailExpon{a}$, $\TailEvent{a}$, and $\lambda(a)$ are given by Equations \eqref{eq: def max queue length}, \eqref{eq: tail expon}, \eqref{eq: tailevent}, and \eqref{subeq: def lambda} respectively, then
\begin{align}\label{eq: 0<a< a star liminf}
\liminf_{N\to\infty}N^{\TailExpon{a}}(\log N)^{\frac{\lambda(a)}{1-\lambda(a)}\frac{\sigma^2}{2\sigma_A^2}}\probability{\MaxQueue{}{}>\TailEvent{a}}>0,    
\end{align}
and
\begin{align}\label{eq: 0<a< a star limsup}
\limsup_{N\to\infty}N^{\TailExpon{a}}(\log N)^{\frac{\lambda(a)}{1-\lambda(a)}\frac{\sigma^2}{2\sigma_A^2}}\probability{\MaxQueue{}{}>\TailEvent{a}}<\infty.    
\end{align}
\end{theorem}

We give a proof of this result in Section \ref{subsec: 0<a< a star}.
As already suggested in Theorem \ref{thm: log asymp}, for the case $0<a<a^{\star}$ we observe more irregular behavior, which manifests itself already in the values of $\TailExpon{a}$. In Theorem \ref{thm: exact asymp 0<a< a star}, we observe that the second term is not a constant, as was the case for the values $a>a^{\star}$ and $a=a^{\star}$, but is $(\log N)^{\frac{\lambda(a)}{1-\lambda(a)}\frac{\sigma^2}{2\sigma_A^2}}$. To obtain heuristic insights, we argue that 
\begin{align}\label{eq: 0<a<a^* sup WA}
    \probability*{\sup_{s>0}(W_A(s)-\lambda(a)\beta s)>\lambda(a)\TailEvent{a}+r_N}=N^{-\TailExpon{a}}(\log N)^{-\frac{\lambda(a)}{1-\lambda(a)}\frac{\sigma^2}{2\sigma_A^2}},
\end{align}
with $r_N=\frac{\sigma  \sqrt{2 a \beta +\sigma ^2}}{4 \beta }\log\log N$. Furthermore, we have for all $k$ that 
\begin{align}
  \probability*{\max_{i\leq N}W_i\big(T_N(a,k)\big)-(1-\lambda(a))\beta T_N(a,k)>(1-\lambda(a))\TailEvent{a}-r_N}=\Omega(1),
\end{align}
where $z_N=\Omega(1)$ means that $\liminf_{N\to\infty}z_N>0$. Combining these two results we see that 
\begin{multline}
    \probability*{\MaxQueue{}{}{}>\TailEvent{a}}\\
    \geq \probability*{\sup_{s>0}(W_A(s)-\lambda(a)\beta s)>\lambda(a)\TailEvent{a}+r_N,\max_{i\leq N}W_i(\tau_N)-(1-\lambda(a))\beta \tau_N>(1-\lambda(a))\TailEvent{a}-r_N},
\end{multline}
where $\tau_N=\inf\{t>0:W_A(t)-\lambda(a)\beta t>\lambda(a)\TailEvent{a}+r_N\}$. We show later on that $\tau_N$ conditioned being finite, has the form of $T_N(a,K)$ with $K$ being a random variable. Because 
\begin{align}
\begin{split}\label{eq: heuristic 0<a<a^*}
    &\probability*{\sup_{s>0}(W_A(s)-\lambda(a)\beta s)>\lambda(a)\TailEvent{a}+r_N,\max_{i\leq N}W_i(\tau_N)-(1-\lambda(a))\beta \tau_N>(1-\lambda(a))\TailEvent{a}-r_N}\\
   &\quad =  \probability*{\sup_{s>0}(W_A(s)-\lambda(a)\beta s)>\lambda(a)\TailEvent{a}+r_N}\\
   &\quad\quad\cdot\probability*{\max_{i\leq N}W_i(\tau_N)-(1-\lambda(a))\beta \tau_N>(1-\lambda(a))\TailEvent{a}-r_N\bigg|\tau_N<\infty},
   \end{split}
\end{align}
we retrieve \eqref{eq: 0<a< a star liminf} after combining the results from \eqref{eq: 0<a<a^* sup WA}--\eqref{eq: heuristic 0<a<a^*}. Thus, it turns out that for $0<a<a^{\star}$, $r_N$ plays a key role. As explained in Section \ref{subsec: a = a star}, in the case $0<a<a^{\star}$, $\{W_A(t)-\lambda(a)\beta t,t>0\}$ dominates, which explains why the tail asymptotics of the maximum queue length $\MaxQueue{}{}$ are the same as the tail asymptotics of $\sup_{s>0}(W_A(s)-\lambda(a)\beta s)$, and the behavior of $\max_{i\leq N}W_i\big(T_N(a,k)\big)-(1-\lambda(a))\beta T_N(a,k)$ is typical.

The main approach of proving the lower and upper bounds in \eqref{eq: 0<a< a star liminf} and \eqref{eq: 0<a< a star limsup}, as well as the limits in \eqref{eq: a>a star thm} and \eqref{eq: a=a star thm}, is by analyzing lower and upper bounds on the tail probability of the steady-state maximum queue length $\probability{\MaxQueue{}{}{}>\TailEvent{a}}$. These bounds are derived by utilizing the union bound, Bonferroni's inequality, and a careful construction of hitting times. These hitting times are needed to estimate the time where the supremum most likely hits the desired level, and to adequately separate the independent part $W_i$ and the dependent part $W_A$ from each other. We also rely on some existing asymptotic estimates in the literature from extreme value theory, and on \cite{debicki2020exact}, that investigates the case $N=2$. Finally, we develop a number of auxiliary technical estimates
related to the asymptotic behavior of convolutions of normally and exponentially distributed random variables.

These techniques, when put together, are effective in the case $a=a^{\star}$ and $a>a^{\star}$ in order to obtain exact asymptotics. 
In the case $0<a<a^{\star}$, we are able to improve upon Theorem \ref{thm: log asymp} and characterize the asymptotic behavior of
$\probability{\MaxQueue{}{}{}>\TailEvent{a}}$ up to a constant. To derive precise asymptotics in this case seems beyond the scope of
techniques developed in this paper.

\section{Proof of the logarithmic asymptotics}\label{sec: log asymp}
In this section, we give a proof of Theorem \ref{thm: log asymp}, establishing logarithmic asymptotics for the maximum queue length. Our approach is to derive logarithmic lower and upper bounds of the maximum queue length by using the heuristic idea given in \eqref{eq: approximate log asymptotics}, and show that they coincide. These bounds are presented in Lemmas \ref{lem: liminf log asymp} and \ref{lem: limsup log asymp} below. 

\begin{lemma}\label{lem: liminf log asymp}
Let $a>0$, $(W_i,i\leq N)$ be i.i.d.\ Brownian motions with standard deviation $\sigma$, $W_A$ be a Brownian motion with standard deviation $\sigma_A$, for all $i$, $W_i$ and $W_A$ are mutually independent, and $\MaxQueue{}{}$, $\TailExpon{a}$, and $\TailEvent{a}$ are given by Equations \eqref{eq: def max queue length}, \eqref{eq: tail expon}, and \eqref{eq: tailevent}, respectively, then
\begin{align}\label{eq: liminf log}
     \liminf_{N\to\infty}\frac{\log(\probability{\MaxQueue{}{}>\TailEvent{a}})}{\log N}\geq -\TailExpon{a}.
\end{align}
\end{lemma}
\begin{proof}
Recall that $\lambda(a)=1-\sigma/\sqrt{2a\beta+\sigma^2}$ and $T_N(a)=\TailEvent{a}/\beta$. By choosing $s=\TailEvent{a}/\beta$ and splitting $-\beta s$ into two terms, observe that 
\begin{align}
&\probability*{\max_{i\leq N}\sup_{s>0}\left(\Process{i}{s}\right)>\TailEvent{a}}\\
&\quad\geq \probability*{\max_{i\leq N}W_i\big(T_N(a)\big)-(1-\lambda(a))\beta T_N(a)>(1-\lambda(a))\TailEvent{a},W_A\big(T_N(a)\big)-\lambda(a)\beta T_N(a)>\lambda(a)\TailEvent{a}}\nonumber\\
&\quad=\probability*{\max_{i\leq N}W_i\big(T_N(a)\big)>2(1-\lambda(a))\TailEvent{a}}\probability*{W_A\big(T_N(a)\big)>2\lambda(a)\TailEvent{a}}.\label{subeq: indep log asymp}
\end{align}
The expression in \eqref{subeq: indep log asymp} is due to the fact that for all $i$, $W_i$ and $W_A$ are independent. We now analyze the two probabilities in \eqref{subeq: indep log asymp} separately. Since $W_i$ and $W_j$ are i.i.d.\ for all $i$ and $j$, for the first probability in \eqref{subeq: indep log asymp} we get from  Bonferroni's inequality that
\begin{multline}
    \probability*{\max_{i\leq N}W_i\big(T_N(a)\big)>2(1-\lambda(a))\TailEvent{a}}\label{eq: bonf lower bound log}\\
    \geq N\probability*{W_i\big(T_N(a)\big)>2(1-\lambda(a))\TailEvent{a}}
    -\binom{N}{2}\probability*{W_i\big(T_N(a)\big)>2(1-\lambda(a))\TailEvent{a}}^2.
    \end{multline}
Furthermore, it is easy to see that 
\begin{align}\label{eq: independent part most likely event}
\probability*{\sup_{s>0}(W_i(s)-(1-\lambda(a))\beta s)>(1-\lambda(a))\TailEvent{a}}=\frac{1}{N}
\end{align}
and that
$$\probability{W_i\big(T_N(a)\big)>2(1-\lambda(a))\TailEvent{a}}\leq \probability*{\sup_{s>0}(W_i(s)-(1-\lambda(a))\beta s)>(1-\lambda(a))\TailEvent{a}},$$ and therefore we bound the second term in \eqref{eq: bonf lower bound log} as 
\begin{align*}
    \binom{N}{2}\probability*{W_i\big(T_N(a)\big)>2(1-\lambda(a))\TailEvent{a}}^2\leq& \frac{N^2}{2}\probability*{\sup_{s>0}(W_i(s)-(1-\lambda(a))\beta s)>(1-\lambda(a))\TailEvent{a}}\nonumber\\
    &\cdot\probability*{W_i\big(T_N(a)\big)>2(1-\lambda(a))\TailEvent{a}}\nonumber\\
    =&\frac{N}{2}\probability*{W_i\big(T_N(a)\big)>2(1-\lambda(a))\TailEvent{a}}.
\end{align*}
Thus the lower bound given in \eqref{eq: bonf lower bound log} can be further bounded to $$\probability*{\max_{i\leq N}W_i\big(T_N(a)\big)>2(1-\lambda(a))\TailEvent{a}}\geq \frac{N}{2}\probability{W_i\big(T_N(a)\big)>2(1-\lambda(a))\TailEvent{a}}.$$ As we aim to derive logarithmic asymptotics, we do so for the derived lower bound, now it is easy to see that 
\begin{align*}
  \log\bigg(\frac{N}{2} \probability*{W_i\big(T_N(a)\big)>2(1-\lambda(a))\TailEvent{a}}\bigg)\sim \log N+\log\bigg( \probability*{W_i\big(T_N(a)\big)>2(1-\lambda(a))\TailEvent{a}}\bigg),
\end{align*}
as $N\to\infty$, with $f(x)\sim g(x)$ as $x\to\infty$ meaning that $\lim_{x\to\infty}f(x)/g(x)=1$. In addition, recall that for a normally distributed random variable $X$ with standard deviation $\sigma$, $\log(\probability*{X>x})\sim -x^2/(2\sigma^2)$, as $x\to\infty$. Thus, we get that
\begin{align*}
    \log\bigg( \probability*{W_i\big(T_N(a)\big)>2(1-\lambda(a))\TailEvent{a}}\bigg)
    \sim -\frac{(2(1-\lambda(a))\TailEvent{a})^2}{2\sigma^2T_N(a)}=-\log N,
\end{align*}
as $N\to\infty$, following the definitions of $\lambda(a)$, $\TailEvent{a}$, and $T_N(a)$. Concluding,
\begin{align}\label{eq: liminf log Wi}
\liminf_{N\to\infty}\frac{\log\bigg(\probability*{\max_{i\leq N}W_i\big(T_N(a)\big)-(1-\lambda(a))\beta  T_N(a)>(1-\lambda(a))\TailEvent{a}}\bigg)}{\log N}\geq 0.
\end{align}

For the second probability in \eqref{subeq: indep log asymp} the logarithmic asymptotics can be easily computed, since $W_A\big(\TailEvent{a}\big)$ is normally distributed, and we obtain that
\begin{align}\label{eq: liminf log WA}
\frac{\log\bigg(\probability*{W_A\big(T_N(a)\big)>2\lambda(a)\TailEvent{a}}\bigg)}{\log N}\LimitN -\frac{2a\beta+2\sigma^2-2\sigma\sqrt{2a\beta+\sigma^2}}{\sigma_A^2}.
\end{align}
Thus, after combining these two results in \eqref{eq: liminf log Wi} and \eqref{eq: liminf log WA} with Equation \eqref{subeq: indep log asymp}, we have that,
\begin{align}\label{eq: liminf log 1}
\liminf_{N\to\infty}\frac{\log\left(\probability*{\max_{i\leq N}\sup_{s>0}\left(\Process{i}{s}\right)>\TailEvent{a}}\right)}{\log N}\geq -\frac{2a\beta+2\sigma^2-2\sigma\sqrt{2a\beta+\sigma^2}}{\sigma_A^2},
\end{align}
irrespective for the choice of $a$. Now, observe that for $a>0$, $$\frac{2a\beta+2\sigma^2-2\sigma\sqrt{2a\beta+\sigma^2}}{\sigma_A^2}\geq\frac{2a\beta-\sigma_A^2}{\sigma^2+\sigma_A^2},$$ with equality for $a=a^{\star}$. This means that only for $0<a\leq a^{\star}$,
the lower bound in \eqref{eq: liminf log 1} is sharp enough. For $a>a^{\star}$, we apply the inequality in \eqref{eq: heuristic ineq memoryless} to obtain for all $c>0$ that
\begin{multline}\label{eq: liminf log 2}
\probability*{\max_{i\leq N}\sup_{s>0}\left(\Process{i}{s}\right)>\TailEvent{a^{\star}+c}}\\
\quad\geq\probability*{\max_{i\leq N}\sup_{s>0}\left(\Process{i}{s}\right)>\TailEvent{a^{\star}}}\exp\left(-\frac{2\beta c\log N}{\sigma^2+\sigma_A^2}\right).
\end{multline}
Combining this result with the inequality in \eqref{eq: liminf log 1}, we get that for all $c>0$,
\begin{align*}
    \liminf_{N\to\infty}\frac{\log\left(\probability*{\max_{i\leq N}\sup_{s>0}\left(\Process{i}{s}\right)>\TailEvent{a^{\star}+c}}\right)}{\log N}\geq -\TailExpon{a^{\star}}-\frac{2\beta c}{\sigma^2+\sigma_A^2}=-\TailExpon{a^{\star}+c}.
\end{align*}
Combining the lower bounds in \eqref{eq: liminf log 1} and \eqref{eq: liminf log 2} gives the lower bound in \eqref{eq: liminf log}. 
\end{proof}
\begin{lemma}\label{lem: limsup log asymp}
Let $a>0$, $(W_i,i\leq N)$ be i.i.d.\ Brownian motions with standard deviation $\sigma$, $W_A$ be a Brownian motion with standard deviation $\sigma_A$, for all $i$, $W_i$ and $W_A$ are mutually independent, and $\MaxQueue{}{}$, $\TailExpon{a}$, and $\TailEvent{a}$ are given by Equations \eqref{eq: def max queue length}, \eqref{eq: tail expon}, and \eqref{eq: tailevent}, respectively, then
\begin{align}\label{eq: limsup log}
     \limsup_{N\to\infty}\frac{\log(\probability{\MaxQueue{}{}>\TailEvent{a}})}{\log N}\leq -\TailExpon{a}.
\end{align}
\end{lemma}
\begin{proof}
We have by the union bound in \eqref{eq: union bound} that
\begin{align}\label{eq: limsup log 1}
    \limsup_{N\to\infty}\frac{\log(\probability{\MaxQueue{}{}>\TailEvent{a}})}{\log N}\leq -\frac{2a\beta-\sigma_A^2}{\sigma^2+\sigma_A^2}.
\end{align}
This upper bound implies the upper bound given in \eqref{eq: limsup log} for $a\geq a^{\star}$. Turning to the case $0<a<a^{\star}$, we can bound the tail probability of the maximum queue length by using sub-additivity, the union bound, and by integrating over possible values of $\sup_{s>0}(W_A(s)-\lambda(a)\beta s)$, and we obtain that
\begin{align}\label{eq: limsup log 2}
&\probability*{\MaxQueue{}{}>\TailEvent{a}}\\
&\quad \leq  \probability*{\max_{i\leq N}\sup_{s>0}\left(W_i(s)-(1-\lambda(a))\beta s\right)+\sup_{s>0}\left(W_A(s)-\lambda(a)\beta s\right)>\TailEvent{a}}\nonumber\\
&\quad \leq \int_{0}^{\lambda(a)(\frac{\sigma^2}{2\beta}+a)}\frac{2\lambda(a)\beta}{\sigma_A^2} N \log N\probability*{\sup_{s>0}\left(W_i(s)-(1-\lambda(a))\beta s\right)>f_N(a)-y\log N}\exp\left(-\frac{2\lambda(a)\beta y\log N}{\sigma_A^2}\right)dy  \nonumber\\
&\quad\quad +\probability*{\sup_{s>0}\left(W_A(s)-\lambda(a)\beta s\right)>\lambda(a)\TailEvent{a}}\nonumber\\
&\quad =\int_{0}^{\lambda(a)(\frac{\sigma^2}{2\beta}+a)}\frac{2\lambda(a)\beta}{\sigma_A^2}N\log N\exp\left(-\frac{2(1-\lambda(a))\beta}{\sigma^2}\left(\TailEvent{a}-y\log N\right)-\frac{2\lambda(a)\beta y\log N}{\sigma_A^2}\right)dy\nonumber\\
&\quad\quad +\probability*{\sup_{s>0}\left(W_A(s)-\lambda(a)\beta s\right)>\lambda(a)\TailEvent{a}}.\nonumber
\end{align}
Because the function $\exp\left(-\frac{2(1-\lambda(a))\beta}{\sigma^2}\left(\TailEvent{a}-y\log N\right)-\frac{2\lambda(a)\beta y\log N}{\sigma_A^2}\right)$ with $y\in[0,\lambda(a)(\frac{\sigma^2}{2\beta}+a)]$ is maximized when $y=\lambda(a)(\frac{\sigma^2}{2\beta}+a)$ and equals $N^{-\frac{2a\beta+2\sigma^2-2\sigma\sqrt{2a\beta+\sigma^2}}{\sigma_A^2}-1}$ we get that
\begin{align}
&\limsup_{N\to\infty}\frac{\log\left(\int_{0}^{\lambda(a)(\frac{\sigma^2}{2\beta}+a)}\frac{2\lambda(a)\beta}{\sigma_A^2}\log N\cdot N\exp\left(-\frac{2(1-\lambda(a))\beta}{\sigma^2}\left(\TailEvent{a}-y\log N\right)-\frac{2\lambda(a)\beta y\log N}{\sigma_A^2}\right)dy\right)}{\log N}\nonumber\\
&\quad =1+\limsup_{N\to\infty}\frac{\log\left(\int_{0}^{\lambda(a)(\frac{\sigma^2}{2\beta}+a)}\exp\left(-\frac{2(1-\lambda(a))\beta}{\sigma^2}\left(\TailEvent{a}-y\log N\right)-\frac{2\lambda(a)\beta y\log N}{\sigma_A^2}\right)dy\right)}{\log N}\nonumber\\
&\quad \leq-\frac{2a\beta+2\sigma^2-2\sigma\sqrt{2a\beta+\sigma^2}}{\sigma_A^2}.
\end{align}
Now we have found a logarithmic upper bound for the integral in \eqref{eq: limsup log 2}, we are left with the expression $\probability*{\sup_{s>0}\left(W_A(s)-\lambda(a)\beta s\right)>\lambda(a)\TailEvent{a}}$ in \eqref{eq: limsup log 2}. For this expression holds that
\begin{align*}
\probability*{\sup_{s>0}\left(W_A(s)-\lambda(a)\beta s\right)>\lambda(a)\TailEvent{a}}=N^{-\frac{2a\beta+2\sigma^2-2\sigma\sqrt{2a\beta+\sigma^2}}{\sigma_A^2}}.    
\end{align*}
Combining the upper bounds in \eqref{eq: limsup log 1} and \eqref{eq: limsup log 2} gives the logarithmic upper bound on the maximum queue length in \eqref{eq: limsup log}. 
\end{proof}

\section{Useful lemmas}\label{sec: useful lemmas}
In the previous section, we have given a proof of the logarithmic asymptotics for the maximum queue length $\MaxQueue{}{}$. In order to be able to prove sharper results on the tail asymptotics, we need some auxiliary results; the goal of this section is to derive these. We begin by giving an overview of the results in this section.

First of all, observe that 
$$
    \sup_{s>T}(W(s)-\beta s)=W(T)-\beta T+\sup_{s>0}(\hat{W}(s)-\beta s),
$$
where $\{\hat{W}(t),t>0\}$ is an independent copy of $\{W(t),t>0\}$. From this, it follows that if we take the supremum of a Brownian motion starting at a positive time, this is in distribution the same as adding a normally distributed random variable to an exponentially distributed random variable. The tail asymptotics of this convolution equal the tail asymptotics of the normally distributed part, the exponentially distributed part, or a more complicated mixture of the two, depending on the starting time $T$, the standard deviation of $W(s)$ and the drift $\beta$. In Lemma \ref{lem: convolution normal exponential}, these three cases are studied in more detail.

Second, our main strategy to investigate the tail asymptotics involves the use of hitting times. Observe that we have a maximum of $N$ mutually dependent random variables. Based on the results in Section \ref{sec: log asymp}, we are able to make an educated guess where the supremum is attained. Following the proof of Lemma \ref{lem: liminf log asymp}, we see that
$$
    \probability*{\max_{i\leq N}\sup_{s>0}\left(\Process{i}{s}\right)>\TailEvent{a}}\approx \probability*{\max_{i\leq N}\left(\Process{i}{T_N(a)}\right)>\TailEvent{a}}.
$$
So the expected hitting time, conditioned on being finite, is approximately $T_N(a)$. Next, observe that for $0<a\leq a^{\star}$,
\begin{align}\label{eq: useful lemmas eq 1}
     \probability*{\max_{i\leq N}\sup_{s>0}\left(W_i(s)-(1-\lambda(a))\beta s\right)>(1-\lambda(a))\TailEvent{a}}=&1-\left(1-\exp\left(-\frac{2(1-\lambda(a))\beta}{\sigma^2}(1-\lambda(a))\TailEvent{a}\right)\right)^N\nonumber\\
     =&\Omega(1),
\end{align}
and
\begin{align}\label{eq: useful lemmas eq 2}
       \probability*{\sup_{s>0}\left(W_A(s)-\lambda(a)\beta s\right)>\lambda(a)\TailEvent{a}}=\exp\left(-\frac{2\lambda(a)\beta}{\sigma_A^2}\lambda(a)\TailEvent{a}\right)=N^{-\TailExpon{a}}.
\end{align}
Since the expected conditional hitting time of a level $x$ equals this value $x$ divided by the drift, it is easy to see that in both \eqref{eq: useful lemmas eq 1} and \eqref{eq: useful lemmas eq 2} the expected conditional hitting time equals $T_N(a)$. Thus, this heuristically explains why the processes $\{W_i(t)-(1-\lambda(a))\beta t,t>0\}$ and $\{W_A(t)-\lambda(a)\beta t,t>0\}$ are important. In Definition \ref{def: hitting time} below, we define the hitting time densities of these processes and in Lemma \ref{lem: hitting time density} we show that after proper scaling these densities converge to the densities of normally distributed random variables, corrected with a constant.

Finally, we need to analyze limits of the type
\begin{align}\label{eq: lim integral cond hitting times}
   \lim_{N\to\infty} \int_{-\infty}^{\infty}\probability*{\sup_{s>\tau_N}X_i(s)>y_N\bigg|\tau_N=t}f_{\tau_N}(t)dt,
\end{align}
where $\tau_N$ is a hitting time and $f_{\tau_N}$ its density. In Lemma \ref{lem: conv sequence integrals}, we show that under certain assumptions, we can interchange the integral and the limit, when the integrand is a product of two functions, as is the case in \eqref{eq: lim integral cond hitting times}. The proof of this interchange is similar to the proof of the dominated convergence theorem. 

\begin{lemma}[Convolution of normal and exponential distributions]\label{lem: convolution normal exponential}
Let $X\overset{d}{=}\mathcal{N}(0,1)$ and $E\overset{d}{=} \text{Exp}(1)$ be independent random variables. Let $(\eta_N,N\geq 1)$, $(x_N,N\geq 1)$ be sequences with $\eta_N>0$, $x_N\LimitN \infty$, and $x_N/\eta_N\LimitN \infty$. Furthermore, let $\mu>0$ and $c\in\mathbb{R}$.
Then
\begin{enumerate}
    \item if $\frac{x_N-\mu  \eta_N ^2}{\sqrt{2} \eta_N }\LimitN c$,
    \begin{align}\label{eq: convolution tail 1}
    \probability*{\eta_NX+\frac{1}{\mu}E>x_N}\sim\frac{\eta_N  e^{-\frac{x_N^2}{2 \eta_N ^2}}}{\sqrt{2 \pi } x_N}+\frac{1}{2} e^{\frac{1}{2} \mu  \left(\mu  \eta_N ^2-2 x_N\right)}(1+\text{erf}(c)),
\end{align}
as $N\to\infty$, and $\text{erf}(c)$ the error function,
\item if $\frac{x_N-\mu  \eta_N ^2}{\sqrt{2} \eta_N }\LimitN \infty$,
    \begin{align}\label{eq: convolution tail 2}
    \probability*{\eta_NX+\frac{1}{\mu}E>x_N}\sim\frac{\eta_N  e^{-\frac{x_N^2}{2 \eta_N ^2}}}{\sqrt{2 \pi } x_N}+e^{\frac{1}{2} \mu  \left(\mu  \eta_N ^2-2 x_N\right)},
\end{align}
as $N\to\infty$,
\item and if $\frac{x_N-\mu  \eta_N ^2}{\sqrt{2} \eta_N }\LimitN -\infty$,
    \begin{align}\label{eq: convolution tail 3}
    \probability*{\eta_NX+\frac{1}{\mu}E>x_N}\sim\frac{\eta_N  e^{-\frac{x_N^2}{2 \eta_N ^2}}}{\sqrt{2 \pi } x_N}-\frac{1}{\sqrt{2\pi}}e^{\frac{1}{2} \mu  \left(\mu  \eta_N ^2-2 x_N\right)}\frac{ \eta_N  e^{-\frac{\left(x_N-\mu  \eta_N ^2\right)^2}{2 \eta_N ^2}}}{x_N-\mu 
   \eta_N ^2},
\end{align}
as $N\to\infty$.

\end{enumerate}

\end{lemma}
\begin{proof}
We have
\begin{align}\label{eq: exact expression convolution normal exponential}
\probability*{\eta_NX+\frac{1}{\mu}E>x_N}=\probability*{\eta_NX>x_N}+\int_{-\infty}^{x_N/\eta_N}\probability*{\frac{1}{\mu}E>x_N-\eta_Nz}\frac{e^{-\frac{z^2}{2}}}{\sqrt{2 \pi }}dz.
\end{align}
The first term satisfies
\begin{align*}
    \probability*{\eta_NX>x_N}\sim \frac{\eta_N  e^{-\frac{x_N^2}{2 \eta_N ^2}}}{\sqrt{2 \pi } x_N},
\end{align*}
as $N\to\infty$. Furthermore,
\begin{align*}
  \int_{-\infty}^{x_N/\eta_N}\probability*{\frac{1}{\mu}E>x_N-\eta_Nz}\frac{e^{-\frac{z^2}{2}}}{\sqrt{2 \pi }}dz = \frac{1}{2} e^{\frac{1}{2} \mu  \left(\mu  \eta_N ^2-2 x_N\right)}
   \left(\text{erf}\left(\frac{x_N-\mu  \eta_N ^2}{\sqrt{2} \eta_N }\right)+1\right).
\end{align*}
Observe that $\text{erf}(z)\to 1$, as $z\to\infty$ and $1+\text{erf}(-z)\sim \frac{e^{-z^2}}{\sqrt{\pi } z}$, as $z\to\infty$. The lemma follows.
\end{proof}
\begin{definition}\label{def: hitting time}
For $a>0$, $r\in\mathbb{R}$, and $i\in\{1,2,\ldots,N\}$, we define the random variable $\HTime{i}{a}{-r}$ by
\begin{align*}
    \HTime{i}{a}{-r}:=\inf\{t>0:W_i(t)-(1-\lambda(a))\beta t>(1-\lambda(a))\TailEvent{a}-r\},
\end{align*}
and the function $f_{\HTime{i}{a}{-r}}$ as its density. Furthermore, we write 
\begin{align*}
    \HTime{\wedge}{a}{-r}:=\min_{i\leq N}\HTime{i}{a}{-r}.
\end{align*}
Similarly, we define the random variable $\HTimeA{a}{r}$ by
\begin{align*}
    \HTimeA{a}{r}:=\inf\{t>0:W_A(t)-\lambda(a)\beta t>\lambda(a)\TailEvent{a}+r\},
\end{align*}
and the function $f_{\HTimeA{a}{r}}$ as its density.
\end{definition}
\begin{lemma}[Convergence of hitting time density]\label{lem: hitting time density}
For the density function $f_{\HTime{i}{a}{-r}}$ given in Definition \ref{def: hitting time} and $T_N(a,k)$ given in Equation \eqref{subeq: TN(a,k) def} we have that
\begin{align}\label{eq: limit hitting time density}
    N\sqrt{\log N}f_{\HTime{i}{a}{-r}}\big(T_N(a,k)\big)\LimitN \frac{\beta ^2 \exp \left(\frac{\beta  \left(8 a^2 \beta ^2 r-\beta ^3 k^2 \sigma  \sqrt{2 a \beta +\sigma ^2}+8 a \beta 
   r \sigma ^2+2 r \sigma ^4\right)}{\sigma  \left(2 a \beta +\sigma ^2\right)^{5/2}}\right)}{\sqrt{\pi } \left(2 a \beta
   +\sigma ^2\right)}.
\end{align}
\end{lemma}
\begin{proof}
The density $f_{\HTime{i}{a}{-r}}(t)$ satisfies
\begin{align*}
    f_{\HTime{i}{a}{-r}}(t)=\frac{(1-\lambda(a))\TailEvent{a}-r}{\sqrt{2\pi}\sigma t^{3/2}}\exp\bigg(-\frac{((1-\lambda(a))\TailEvent{a}-r+(1-\lambda(a))\beta t)^2}{2\sigma^2 t}\bigg); 
\end{align*}
cf.\ \cite[Eq.\ 2.0.2, p.\ 301]{borodin2015handbook}. From this, the limit in \eqref{eq: limit hitting time density} follows.
\end{proof}
\begin{corollary}
For the density function $f_{\HTime{i}{a}{-r}}$ given in Definition \ref{def: hitting time} and $T_N(a,k)$ given in Equation \eqref{subeq: TN(a,k) def} we have that
\begin{align}
    \lim_{N\to\infty}\int_{-\infty}^{\infty}N\sqrt{\log N}f_{\HTime{i}{a}{-r}}\big(T_N(a,k)\big)dk=\int_{-\infty}^{\infty} \lim_{N\to\infty}N\sqrt{\log N}f_{\HTime{i}{a}{-r}}\big(T_N(a,k)\big)dk.
\end{align}
\end{corollary}
\begin{proof}
Observe that for $N$ large enough such that $(1-\lambda(a))\TailEvent{a}-r>0$,
\begin{align*}
    \int_{-\infty}^{\infty}N\sqrt{\log N}f_{\HTime{i}{a}{-r}}\big(T_N(a,k)\big)dk=&N\probability*{\sup_{s>0}(W_i(s)-(1-\lambda(a))\beta s)>(1-\lambda(a))\TailEvent{a}-r}\nonumber\\
    =&\exp\left(\frac{2(1-\lambda(a))\beta r}{\sigma^2}\right),
\end{align*}
and
\begin{align*}
     \int_{-\infty}^{\infty}\frac{\beta ^2 \exp \left(\frac{\beta  \left(8 a^2 \beta ^2 r-\beta ^3 k^2 \sigma  \sqrt{2 a \beta +\sigma ^2}+8 a \beta 
   r \sigma ^2+2 r \sigma ^4\right)}{\sigma  \left(2 a \beta +\sigma ^2\right)^{5/2}}\right)}{\sqrt{\pi } \left(2 a \beta
   +\sigma ^2\right)}dk=\exp\left(\frac{2(1-\lambda(a))\beta r}{\sigma^2}\right).
\end{align*}
\end{proof}

\begin{lemma}[Convergence of integrals of sequences of functions]\label{lem: conv sequence integrals}
Assume we have sequences of positive integrable functions $v_N(x)$ and $w_N(x)$ that satisfy the following:
\begin{itemize}
    \item $v_N(x)\LimitN v(x)$, 
    \item $\int_{\mathbb{R}} v_N(x)dx\LimitN\int_{\mathbb{R}} v(x)dx$,
    \item $w_N(x)\LimitN w(x)$,
    \item There exists a constant $c>0$ such that $w_N(x)<c$ for all $x$ and $N$.
\end{itemize}
Then 
\begin{align}
    \int_{\mathbb{R}} v_N(x)w_N(x)dx\LimitN \int_{\mathbb{R}} v(x)w(x)dx.
\end{align}
\end{lemma}
\begin{proof}
First of all, by using Fatou's lemma we obtain that
\begin{align*}
    \liminf_{N\to\infty}\int_{\mathbb{R}} v_N(x)w_N(x)dx\geq \int_{\mathbb{R}} v(x)w(x)dx.
\end{align*}
Furthermore, observe that $v_N(x)c-v_N(x)w_N(x)>0$ for all $x$ and $N$. Now, from Fatou's lemma it follows that
\begin{align*}
    \liminf_{N\to\infty}\int_{\mathbb{R}} v_N(x)c-v_N(x)w_N(x)dx\geq \int_{\mathbb{R}} v(x)c-v(x)w(x)dx.
\end{align*}
Because $\int_{\mathbb{R}} v_N(x)cdx\LimitN\int_{\mathbb{R}} v(x)cdx$, we get that
\begin{align*}
    \limsup_{N\to\infty}\int_{\mathbb{R}} v_N(x)w_N(x)dx\leq \int_{\mathbb{R}} v(x)w(x)dx.
\end{align*}
The lemma follows.
\end{proof}

\section{Proofs of the sharper asymptotics}\label{sec: exact asymptotics}
In this section, we prove sharper asymptotics of the tail behavior of $\probability{\MaxQueue{}{}>\TailEvent{a}}$. Recall the definition of $\HTime{i}{a}{-r}$ and $\HTimeA{a}{r}$
given in Definition \ref{def: hitting time}, and 
observe that 
\begin{align}\label{eq: repres max Q hitting times}
    \probability{\MaxQueue{}{}>\TailEvent{a}}
    =\probability{\MaxQueue{\HTime{i}{a}{-r}\wedge \HTimeA{a}{r}}{}\mathbbm{1}(\HTime{\wedge}{a}{-r}\wedge \HTimeA{a}{r}<\infty)>\TailEvent{a}}.
\end{align}
This equation is valid, because for $0<t<\HTime{\wedge}{a}{-r}\wedge \HTimeA{a}{r}$, we see that $W_i(t)-(1-\lambda(a))\beta t<(1-\lambda(a))\TailEvent{a}-r$ and $W_A(t)-\lambda(a)\beta t<\lambda(a)\TailEvent{a}+r$. Thus, $W_i(t)+W_A(t)-\beta t<\TailEvent{a}$. Now, using (\ref{eq: repres max Q hitting times}), we obtain  lower and upper bounds of the form
\begin{align}\label{eq: upper lower bounds repr max Q hitting times}
\max\bigg(&\mathbb{P}\bigg(\MaxQueue{\HTime{i}{a}{-r}}{}\mathbbm{1}(\HTime{\wedge}{a}{-r}<\infty)>\TailEvent{a}\bigg),\mathbb{P}\bigg(\MaxQueue{\HTimeA{a}{r}}{}\mathbbm{1}( \HTimeA{a}{r}<\infty)>\TailEvent{a}\bigg)\bigg)\nonumber\\
\leq&\probability{\MaxQueue{}{}>\TailEvent{a}}\nonumber\\
   \leq&\probability*{\MaxQueue{\HTime{i}{a}{-r}}{}\mathbbm{1}(\HTime{\wedge}{a}{-r}<\infty)>\TailEvent{a}}
   +\probability*{\MaxQueue{\HTimeA{a}{r}}{}\mathbbm{1}( \HTimeA{a}{r}<\infty)>\TailEvent{a}},
\end{align}
which we can exploit. 
Other important inequalities that we use are the union bound and Bonferroni's inequality. In the case of identically distributed random variables $X_i$, these bounds simplify to
\begin{align*}
    N\probability{X_i>x}-\binom{N}{2}\probability{\min(X_i,X_j)>x}\leq \probability{\max_{i\leq N}X_i>x}\leq N\probability{X_i>x},
\end{align*}
which is the case for our problem. D\k{e}bicki et al.\ \cite{debicki2020exact} have derived the tail asymptotics of $\min(\BrQueueLength{i}{}{},\BrQueueLength{j}{}{})$. In Lemma \ref{lem: liminf a > a star} we show how we use \cite[Th. 2.3]{debicki2020exact} on the tails of $\min(\BrQueueLength{i}{}{},\BrQueueLength{j}{}{})$ together with Bonferroni's inequality such that these are applicable in our proof of the case $a>a^{\star}$.

Now that we can write upper and lower bounds in which hitting times play a role, we condition on the hitting times and get sequences of the form as given in \eqref{eq: lim integral cond hitting times}. By using Fatou's lemma we know that 
\begin{align*}
     \liminf_{N\to\infty}\int_{-\infty}^{\infty}\probability*{\sup_{s>\tau_N}X_i(s)>y_N\Big|\tau_N=t}f_{\tau_N}(t)dt\geq \int_{-\infty}^{\infty}\liminf_{N\to\infty}\probability*{\sup_{s>\tau_N}X_i(s)>y_N\Big|\tau_N=t}f_{\tau_N}(t)dt,
\end{align*}
and by using Lemma \ref{lem: conv sequence integrals}, we obtain that 
\begin{align*}
     \lim_{N\to\infty}\int_{-\infty}^{\infty}\probability*{\sup_{s>\tau_N}X_i(s)>y_N\Big|\tau_N=t}f_{\tau_N}(t)dt= \int_{-\infty}^{\infty}\lim_{N\to\infty}\probability*{\sup_{s>\tau_N}X_i(s)>y_N\Big|\tau_N=t}f_{\tau_N}(t)dt.
\end{align*}
To obtain limits of the form as given in \eqref{eq: lim integral cond hitting times} we use Lemmas \ref{lem: convolution normal exponential} and \ref{lem: hitting time density}.

\subsection{The case $a>a^{\star}$}\label{subsec: a > a star}
In this section, we prove Theorem \ref{thm: exact asymp a > a star} on exact asymptotics of the maximum queue length when $a>a^{\star}$. As is stated in \eqref{eq: a>a star thm}, $\probability{\MaxQueue{}{}>\TailEvent{a}}\sim N^{-\TailExpon{a}}$, as $N\to\infty$, when $a>a^{\star}$. Since the union bound in \eqref{eq: union bound} gives us that $N^{\TailExpon{a}}\probability{\MaxQueue{}{}>\TailEvent{a}}\leq 1$,
we only need to show that 
\begin{align*}
    \liminf_{N\to\infty}     N^{\TailExpon{a}}\probability{\MaxQueue{}{}>\TailEvent{a}}\geq 1.
\end{align*}
In order to prove the $\liminf$, we first observe that $\MaxQueue{}{}>\MaxQueue{\HTime{i}{a^{\star}}{0}}{}\mathbbm{1}(\HTime{\wedge}{a^{\star}}{0}<\infty)$, and we know by using Bonferroni's inequality that
\begin{align}\label{eq: a > a star bonf}
&\probability{\MaxQueue{\HTime{i}{a^{\star}}{0}}{}\mathbbm{1}(\HTime{\wedge}{a^{\star}}{0}<\infty)>\TailEvent{a}}\nonumber\\
&\quad\geq N\probability*{\BrQueueLength{i}{\HTime{i}{a^{\star}}{0}}{}\mathbbm{1}(\HTime{i}{a^{\star}}{0}<\infty)>\TailEvent{a}}
\nonumber\\
   &\quad\quad-\binom{N}{2} \probability*{\min(\BrQueueLength{i}{\HTime{i}{a^{\star}}{0}}{}\mathbbm{1}(\HTime{i}{a^{\star}}{0}<\infty),\BrQueueLength{j}{\HTime{i}{a^{\star}}{0}}{}\mathbbm{1}(\HTime{j}{a^{\star}}{0}<\infty))>\TailEvent{a}},
\end{align}
where $\HTime{i}{a^{\star}}{0}$ and $\HTime{j}{a^{\star}}{0}$ are hitting times defined in Lemma \ref{lem: hitting time density}. In Lemma \ref{lem: liminf a > a star}, we show that the first term is leading, and the second order term is of smaller order. 
In order to prove this, we first give a convenient upper bound for
$$
  \CondProb{k}{l}*{\min(\BrQueueLength{i}{\HTime{i}{a^{\star}}{0}}{},\BrQueueLength{j}{\HTime{j}{a^{\star}}{0}}{})>\TailEvent{a}}
$$
in Lemma \ref{lem: upper bound joint prob}, with 
\begin{align}\label{eq: def cond prob}
      \CondProb{k}{l}*{A}=\probability*{A\Big| \HTime{i}{a^{\star}}{0}=T_N(a^{\star},k)<\HTime{j}{a^{\star}}{0}=T_N(a^{\star},l)}.
\end{align}
From now on, let $\hat{W}$ be an independent copy of the Brownian motion $W$, and $\HatBrQueueLength{i}{s}{t}$ an independent copy of $\BrQueueLength{i}{u}{v}$.
\begin{lemma}\label{lem: upper bound joint prob}
Let $a>a^{\star}$, $(W_i,i\leq N)$ be i.i.d.\ Brownian motions with standard deviation $\sigma$, $W_A$ be a Brownian motion with standard deviation $\sigma_A$, for all $i$, $W_i$ and $W_A$ are mutually independent, and $\BrQueueLength{i}{u}{v}$, $\TailExpon{a}$, $\TailEvent{a}$, and $\CondProb{k}{l}{A}$ are given by Equations \eqref{eq: queue length finite time}, \eqref{eq: tail expon}, \eqref{eq: tailevent}, and \eqref{eq: def cond prob} respectively. Furthermore, $\HTime{i}{a^{\star}}{0}$ is given in Definition \ref{def: hitting time} and $\HatBrQueueLength{i}{}{}$ is an independent copy of $\BrQueueLength{i}{}{}$. Then for all $\delta>0$ there exists an $N_{\delta}>0$ such that for all $N\geq N_{\delta}$
\begin{multline*}
    \CondProb{k}{l}*{\min(\BrQueueLength{i}{\HTime{i}{a^{\star}}{0}}{},\BrQueueLength{j}{\HTime{j}{a^{\star}}{0}}{})>\TailEvent{a}}\\
 \leq 4\mathbb{P}^{(k<l)}\bigg((1+\delta)W_A(\HTime{i}{a^{\star}}{0})+\min(\HatBrQueueLength{i}{}{},\HatBrQueueLength{j}{}{})
    >\TailEvent{a}-(1-\lambda(a^{\star}))\TailEvent{a^{\star}}+\lambda(a^{\star})\beta\HTime{i}{a^{\star}}{0}\bigg).
\end{multline*}
\end{lemma}
\begin{proof}
First of all, we have that
\begin{multline}
\CondProb{k}{l}*{\min(\BrQueueLength{i}{\HTime{i}{a^{\star}}{0}}{},\BrQueueLength{j}{\HTime{j}{a^{\star}}{0}}{})>\TailEvent{a}}\\
   \leq  \CondProb{k}{l}*{\BrQueueLength{i}{\HTime{i}{a^{\star}}{0}}{\HTime{j}{a^{\star}}{0}}>\TailEvent{a}}
   +\CondProb{k}{l}*{\min(\BrQueueLength{i}{\HTime{j}{a^{\star}}{0}}{},\BrQueueLength{j}{\HTime{j}{a^{\star}}{0}}{})>\TailEvent{a}},\label{eq: bonf term 1}
\end{multline}
because $\min(\BrQueueLength{i}{\HTime{i}{a^{\star}}{0}}{},\BrQueueLength{j}{\HTime{j}{a^{\star}}{0}}{})<\max(\BrQueueLength{i}{\HTime{i}{a^{\star}}{0}}{\HTime{j}{a^{\star}}{0}},\min(\BrQueueLength{i}{\HTime{j}{a^{\star}}{0}}{},\BrQueueLength{j}{\HTime{j}{a^{\star}}{0}}{}))$ when $\HTime{i}{a^{\star}}{0}<\HTime{j}{a^{\star}}{0}<\infty$. Now, recall from Definition \ref{def: hitting time} that
\begin{multline*}
    \BrQueueLength{i}{\HTime{i}{a^{\star}}{0}}{\HTime{j}{a^{\star}}{0}}\\=\sup_{\HTime{i}{a^{\star}}{0}<s<\HTime{j}{a^{\star}}{0}}(\Process{i}{s})
=(1-\lambda(a^{\star}))\TailEvent{a^{\star}}+W_A(\HTime{i}{a^{\star}}{0})-\lambda(a^{\star})\beta\HTime{i}{a^{\star}}{0}+\HatBrQueueLength{i}{0}{\HTime{j}{a^{\star}}{0}-\HTime{i}{a^{\star}}{0}}.
\end{multline*}
Thus, for the first term in \eqref{eq: bonf term 1} we have
\begin{align}
  &\CondProb{k}{l}*{\BrQueueLength{i}{\HTime{i}{a^{\star}}{0}}{\HTime{j}{a^{\star}}{0}}>\TailEvent{a}}\nonumber\\
  &\quad =\CondProb{k}{l}*{W_A(\HTime{i}{a^{\star}}{0})+\HatBrQueueLength{i}{0}{\HTime{j}{a^{\star}}{0}-\HTime{i}{a^{\star}}{0}}>\TailEvent{a}-(1-\lambda(a^{\star}))\TailEvent{a^{\star}}+\lambda(a^{\star})\beta\HTime{i}{a^{\star}}{0}}\nonumber\\
  &\quad \leq\CondProb{k}{l}*{W_A(\HTime{i}{a^{\star}}{0})+\left|\hat{W}_i(\HTime{j}{a^{\star}}{0}-\HTime{i}{a^{\star}}{0})+\hat{W}_A(\HTime{j}{a^{\star}}{0}-\HTime{i}{a^{\star}}{0})\right|>\TailEvent{a}-(1-\lambda(a^{\star}))\TailEvent{a^{\star}}+\lambda(a^{\star})\beta\HTime{i}{a^{\star}}{0}}\label{eq: bonf term 2}.
\end{align}
For any $x$ and $y$, it holds that $x+|y|\leq \max(x+y,x-y)$. Therefore, by the union bound we can bound the probability in \eqref{eq: bonf term 2} as
\begin{align}
    &\CondProb{k}{l}*{W_A(\HTime{i}{a^{\star}}{0})+\left|\hat{W}_i(\HTime{j}{a^{\star}}{0}-\HTime{i}{a^{\star}}{0})+\hat{W}_A(\HTime{j}{a^{\star}}{0}-\HTime{i}{a^{\star}}{0})\right|>\TailEvent{a}-(1-\lambda(a^{\star}))\TailEvent{a^{\star}}+\lambda(a^{\star})\beta\HTime{i}{a^{\star}}{0}}\nonumber\\
  &\quad \leq  2\CondProb{k}{l}*{W_A(\HTime{i}{a^{\star}}{0})+\hat{W}_i(\HTime{j}{a^{\star}}{0}-\HTime{i}{a^{\star}}{0})+\hat{W}_A(\HTime{j}{a^{\star}}{0}-\HTime{i}{a^{\star}}{0})>\TailEvent{a}-(1-\lambda(a^{\star}))\TailEvent{a^{\star}}+\lambda(a^{\star})\beta\HTime{i}{a^{\star}}{0}}\nonumber\\
  &\quad \leq 2\CondProb{k}{l}*{(1+\delta)W_A(\HTime{i}{a^{\star}}{0})>\TailEvent{a}-(1-\lambda(a^{\star}))\TailEvent{a^{\star}}+\lambda(a^{\star})\beta\HTime{i}{a^{\star}}{0}}\label{subeq: bonferroni bound 1}\\
  &\quad \leq 2\mathbb{P}^{(k<l)}\bigg((1+\delta)W_A(\HTime{i}{a^{\star}}{0})+\min(\HatBrQueueLength{i}{}{},\HatBrQueueLength{j}{}{})
    >\TailEvent{a}-(1-\lambda(a^{\star}))\TailEvent{a^{\star}}+\lambda(a^{\star})\beta\HTime{i}{a^{\star}}{0}\bigg)\label{subeq: bonferroni bound 2},
\end{align}
for $\delta>0$ and $N>N_{\delta}$. The upper bound in \eqref{subeq: bonferroni bound 1} holds since $\HTime{i}{a^{\star}}{0}=\Omega(\log N)$, and $\HTime{j}{a^{\star}}{0}-\HTime{i}{a^{\star}}{0}=O(\sqrt{\log N})$. The upper bound in \eqref{subeq: bonferroni bound 2} holds because we add a positive random variable.
For the second term in \eqref{eq: bonf term 1}, first observe that $\probability{\min(X,Y)>z}=\probability{X>z,Y>z}$. Second, under the assumption that $\HTime{i}{a^{\star}}{0}<\HTime{j}{a^{\star}}{0}<\infty$, we can write 
$$
\BrQueueLength{i}{\HTime{j}{a^{\star}}{0}}{}=(1-\lambda(a^{\star}))\TailEvent{a^{\star}}+W_i(\HTime{j}{a^{\star}}{0}-\HTime{i}{a^{\star}}{0})-(1-\lambda(a^{\star}))\beta(\HTime{j}{a^{\star}}{0}-\HTime{i}{a^{\star}}{0})+W_A(\HTime{j}{a^{\star}}{0})-\lambda(a^{\star})\beta\HTime{j}{a^{\star}}{0}+\HatBrQueueLength{i}{}{}.
$$
Thus, by applying similar techniques as for the analysis of the first term in \eqref{eq: bonf term 1} we obtain that
\begin{align*}
    &\CondProb{k}{l}*{\min(\BrQueueLength{i}{\HTime{j}{a^{\star}}{0}}{},\BrQueueLength{j}{\HTime{j}{a^{\star}}{0}}{})>\TailEvent{a}}\nonumber\\
   &\quad  =\mathbb{P}^{(k<l)}\bigg(W_A(\HTime{j}{a^{\star}}{0})+W_i(\HTime{j}{a^{\star}}{0}-\HTime{i}{a^{\star}}{0})-(1-\lambda(a^{\star}))\beta(\HTime{j}{a^{\star}}{0}-\HTime{i}{a^{\star}}{0})+\HatBrQueueLength{i}{}{}\nonumber\\
    &\quad \hspace{1.5cm}>\TailEvent{a}-(1-\lambda(a^{\star}))\TailEvent{a^{\star}}+\lambda(a^{\star})\beta\HTime{j}{a^{\star}}{0},\nonumber\\
    &\quad \hspace{1cm} W_A(\HTime{j}{a^{\star}}{0})+\HatBrQueueLength{j}{}{}>\TailEvent{a}-(1-\lambda(a^{\star}))\TailEvent{a^{\star}}+\lambda(a^{\star})\beta\HTime{j}{a^{\star}}{0}
    \bigg)\nonumber\\
      &\quad \leq\mathbb{P}^{(k<l)}\bigg(W_A(\HTime{j}{a^{\star}}{0})+W_i(\HTime{j}{a^{\star}}{0}-\HTime{i}{a^{\star}}{0})+\HatBrQueueLength{i}{}{}
    >\TailEvent{a}-(1-\lambda(a^{\star}))\TailEvent{a^{\star}}+\lambda(a^{\star})\beta\HTime{j}{a^{\star}}{0},\nonumber\\
      &\quad \hspace{1cm} W_A(\HTime{j}{a^{\star}}{0})+\HatBrQueueLength{j}{}{}>\TailEvent{a}-(1-\lambda(a^{\star}))\TailEvent{a^{\star}}+\lambda(a^{\star})\beta\HTime{j}{a^{\star}}{0}\bigg)\nonumber\\
      &\quad \leq\mathbb{P}^{(k<l)}\bigg(W_A(\HTime{j}{a^{\star}}{0})+\max(W_i(\HTime{j}{a^{\star}}{0}-\HTime{i}{a^{\star}}{0}),0)+\min(\HatBrQueueLength{i}{}{},\HatBrQueueLength{j}{}{})
    >\TailEvent{a}-(1-\lambda(a^{\star}))\TailEvent{a^{\star}}+\lambda(a^{\star})\beta\HTime{j}{a^{\star}}{0}\bigg)\nonumber\\
      &\quad \leq 2\mathbb{P}^{(k<l)}\bigg(W_A(\HTime{j}{a^{\star}}{0})+W_i(\HTime{j}{a^{\star}}{0}-\HTime{i}{a^{\star}}{0})+\min(\HatBrQueueLength{i}{}{},\HatBrQueueLength{j}{}{})
    >\TailEvent{a}-(1-\lambda(a^{\star}))\TailEvent{a^{\star}}+\lambda(a^{\star})\beta\HTime{j}{a^{\star}}{0}\bigg)\nonumber\\
      &\quad \leq 2\mathbb{P}^{(k<l)}\bigg((1+\delta)W_A(\HTime{j}{a^{\star}}{0})+\min(\HatBrQueueLength{i}{}{},\HatBrQueueLength{j}{}{})>\TailEvent{a}-(1-\lambda(a^{\star}))\TailEvent{a^{\star}}+\lambda(a^{\star})\beta\HTime{i}{a^{\star}}{0}\bigg).
\end{align*}
Combining this bound with the bound in \eqref{subeq: bonferroni bound 2} completes the proof of the lemma.
\end{proof}
\begin{lemma}\label{lem: liminf a > a star}
Let $a>a^{\star}$ $(W_i,i\leq N)$ be i.i.d.\ Brownian motions with standard deviation $\sigma$, $W_A$ be a Brownian motion with standard deviation $\sigma_A$, for all $i$, $W_i$ and $W_A$ are mutually independent, and $\MaxQueue{}{}$, $\TailExpon{a}$, and $\TailEvent{a}$ are given by Equations \eqref{eq: def max queue length}, \eqref{eq: tail expon}, and \eqref{eq: tailevent}, respectively, then 
$$
   \liminf_{N\to\infty} N^{\TailExpon{a}}\probability{\MaxQueue{}{}>\TailEvent{a}}\geq 1.
$$
\end{lemma}

The general idea of the proof of Lemma \ref{lem: liminf a > a star} is to make rigorous that the lower bound on the maximum queue length $\MaxQueue{}{}$ given in \eqref{eq: a > a star bonf} is approximately the same as $N\probability{\BrQueueLength{i}{\HTime{i}{a^{\star}}{0}}{}\mathbbm{1}(\HTime{i}{a^{\star}}{0}<\infty)>\TailEvent{a}}$ when $N$ is large. Thus the last term in \eqref{eq: a > a star bonf} is asymptotically negligible.
We use the result from Lemma \ref{lem: upper bound joint prob} to establish this. Observe now that, following Definition \ref{def: hitting time},
$$\BrQueueLength{i}{\HTime{i}{a^{\star}}{0}}{}=W_i(\HTime{i}{a^{\star}}{0})+W_A(\HTime{i}{a^{\star}}{0})-\beta \HTime{i}{a^{\star}}{0}+\HatBrQueueLength{i}{}{}=(1-\lambda(a^{\star}))\TailEvent{a^{\star}}+W_A(\HTime{i}{a^{\star}}{0})-\lambda(a^{\star})\beta \HTime{i}{a^{\star}}{0}+\HatBrQueueLength{i}{}{}.$$ Furthermore, observe that due to Equation \eqref{eq: independent part most likely event},  $\probability{\HTime{i}{a^{\star}}{0}<\infty}=1/N$. From this it follows that $$ N\probability{\BrQueueLength{i}{\HTime{i}{a^{\star}}{0}}{}\mathbbm{1}(\HTime{i}{a^{\star}}{0}<\infty)>\TailEvent{a}}=\probability{\BrQueueLength{i}{\HTime{i}{a^{\star}}{0}}{}>\TailEvent{a}\mid\HTime{i}{a^{\star}}{0}<\infty}.
$$
Therefore, in order to prove a sharp lower bound on the tail asymptotics of the maximum queue length, we prove by using Fatou's lemma that $$\liminf_{N\to\infty}N^{\TailExpon{a}}\probability{W_A(\HTime{i}{a^{\star}}{0})-\lambda(a^{\star})\beta \HTime{i}{a^{\star}}{0}+\HatBrQueueLength{i}{}{}>\TailEvent{a}-(1-\lambda(a^{\star}))\TailEvent{a^{\star}}\mid\HTime{i}{a^{\star}}{0}<\infty}\geq 1.$$ In order to prove this, we show that $\HatBrQueueLength{i}{}{}$ is most likely to hit a level $g_N(a,x,k)$, and $W_A(\HTime{i}{a^{\star}}{0})-\lambda(a^{\star})\beta \HTime{i}{a^{\star}}{0}$ is most likely to hit the level $\TailEvent{a}-(1-\lambda(a^{\star}))\TailEvent{a^{\star}}-g_N(a,x,k)$.\\

We now turn to a formal proof of Lemma \ref{lem: liminf a > a star}.
\begin{proof}
For abbreviation, we write
$$
    P_{i,j,N}=\probability*{\min(\BrQueueLength{i}{\HTime{i}{a^{\star}}{0}}{}\mathbbm{1}(\HTime{i}{a^{\star}}{0}<\infty),\BrQueueLength{j}{\HTime{j}{a^{\star}}{0}}{}\mathbbm{1}(\HTime{j}{a^{\star}}{0}<\infty))>\TailEvent{a}}.
$$
Thus, the inequality in \eqref{eq: a > a star bonf} simplifies to
\begin{align}\label{eq: bonf inequality}
   \probability{\MaxQueue{\HTime{i}{a^{\star}}{0}}{}\mathbbm{1}(\HTime{\wedge}{a^{\star}}{0}<\infty)>\TailEvent{a}}
    \geq NP_{i,i,N}-\binom{N}{2}P_{i,j,N}.
\end{align}
For abbreviation, we also write
$$
    Q_{i,j,N}(k,l)=
    \probability*{\min(\BrQueueLength{i}{\HTime{i}{a^{\star}}{0}}{},\BrQueueLength{j}{\HTime{j}{a^{\star}}{0}}{})>\TailEvent{a}
    \Big|\HTime{i}{a^{\star}}{0}=T_N(a^{\star},k),\HTime{j}{a^{\star}}{0}=T_N(a^{\star},l)}.
$$
Now, before we analyze \eqref{eq: bonf inequality} in more detail, observe that we can express $\probability{\HTime{i}{a^{\star}}{0}<\infty,\HTime{j}{a^{\star}}{0}<\infty}$ as 
$$
\probability{\HTime{i}{a^{\star}}{0}<\infty,\HTime{j}{a^{\star}}{0}<\infty}=\int_{-\infty}^{\infty}\int_{-\infty}^{\infty}f_{\HTime{i}{a^{\star}}{0}}\big(T_N(a^{\star},k)\big)f_{\HTime{j}{a^{\star}}{0}}\big(T_N(a^{\star},l)\big)\log N dkdl=\frac{1}{N^2}.
$$
Then,
\begin{align*}
    NP_{i,i,N}=&N\int_{-\infty}^{\infty}f_{\HTime{i}{a^{\star}}{0}}\big(T_N(a^{\star},k)\big)\sqrt{\log N}Q_{i,i,N}(k,k)dk\\
   =&\int_{-\infty}^{\infty}\int_{-\infty}^{\infty}f_{\HTime{i}{a^{\star}}{0}}\big(T_N(a^{\star},k)\big)f_{\HTime{j}{a^{\star}}{0}}\big(T_N(a^{\star},l)\big)N^2\log NQ_{i,i,N}(k,k)dkdl.
\end{align*}
Also, observe that $\binom{N}{2}<N^2/2$, and that
\begin{align*}
    \frac{N^2}{2}P_{i,j,N}=\frac{N^2}{2}\int_{-\infty}^{\infty}\int_{-\infty}^{\infty}f_{\HTime{i}{a^{\star}}{0}}\big(T_N(a^{\star},k)\big)f_{\HTime{j}{a^{\star}}{0}}\big(T_N(a^{\star},l)\big)\log NQ_{i,j,N}(k,l)dkdl.
\end{align*}
In conclusion, we can write the inequality in \eqref{eq: bonf inequality} as
\begin{align}
        &\probability{\MaxQueue{\HTime{i}{a^{\star}}{0}}{}\mathbbm{1}(\HTime{\wedge}{a^{\star}}{0}<\infty)>\TailEvent{a}}\label{eq: lower bound a > a star}\\
        &\quad\geq\int_{-\infty}^{\infty}\int_{-\infty}^{\infty}f_{\HTime{i}{a^{\star}}{0}}\big(T_N(a^{\star},k)\big)f_{\HTime{j}{a^{\star}}{0}}\big(T_N(a^{\star},l)\big)N^2\log N\left(Q_{i,i,N}(k,k)-\frac{1}{2}Q_{i,j,N}(k,l)\right)dkdl\label{subeq: lower bound a>a star}\\
        &\quad=\int_{-\infty}^{\infty}\int_{-\infty}^{l}f_{\HTime{i}{a^{\star}}{0}}\big(T_N(a^{\star},k)\big)f_{\HTime{j}{a^{\star}}{0}}\big(T_N(a^{\star},l)\big)N^2\log N\left(Q_{i,i,N}(k,k)-\frac{1}{2}Q_{i,j,N}(k,l)\right)dkdl\nonumber\\
        &\quad\quad+\int_{-\infty}^{\infty}\int_{l}^{\infty}f_{\HTime{i}{a^{\star}}{0}}\big(T_N(a^{\star},k)\big)f_{\HTime{j}{a^{\star}}{0}}\big(T_N(a^{\star},l)\big)N^2\log N\left(Q_{i,i,N}(k,k)-\frac{1}{2}Q_{i,j,N}(k,l)\right)dkdl.\nonumber
\end{align}
Since we want to prove a sharp lower bound on the tail asymptotics of the maximum queue length $\MaxQueue{}{}{}$ we can use the expression in \eqref{subeq: lower bound a>a star}. We want to prove convergence of a lower bound of this integral by using Fatou's lemma. Therefore, we focus on the integrand first and prove convergence for the integrand as $N\to\infty$.
Assume that $k\leq l$, then, following Lemma \ref{lem: upper bound joint prob},
$$
   Q_{i,j,N}(k,l)
  \leq 4\mathbb{P}\bigg((1+\delta)W_A\big(T_N(a^{\star},k)\big)
  +\min(\HatBrQueueLength{i}{}{},\HatBrQueueLength{j}{}{})
  >
  \TailEvent{a}-(1-\lambda(a^{\star}))\TailEvent{a^{\star}}+\lambda(a^{\star})\beta T_N(a^{\star},k)\bigg),
$$
for all $\delta>0$ for $N>N_{\delta}$.
Observe that $Q_{i,i,N}(k,k)-Q_{i,j,N}(k,l)/2>0$. Thus,
$$
    Q_{i,i,N}(k,k)-\frac{1}{2}Q_{i,j,N}(k,l)=\left(Q_{i,i,N}(k,k)-\frac{1}{2}Q_{i,j,N}(k,l)\right)^+.
$$
The density of $W_A\big(T_N(a^{\star},k)\big)$ equals
$$
    \frac{\exp\left(-x^2/(2\sigma_A^2 T_N(a^{\star},k))\right)}{\sqrt{2\pi}\sigma_A\sqrt{T_N(a^{\star},k)}}.
$$
We write $a=a^{\star}+\epsilon$, with $\epsilon>0$. Let 
$$
g_N(a,x,k)
=\TailEvent{a}-(1-\lambda(a^{\star}))\TailEvent{a^{\star}}+\lambda(a^{\star})\beta T_N(a^{\star},k)-\frac{\sigma_A ^2 \left(\sigma ^2+\sigma_A ^2\right)}{\beta  \sigma ^2}\log N-x\sqrt{\log N}.
$$
Observe that 
$$
   g_N(a,x,k)+\frac{\sigma_A ^2 \left(\sigma ^2+\sigma_A ^2\right)}{\beta  \sigma ^2}\log N+x\sqrt{\log N}
   = \TailEvent{a}-(1-\lambda(a^{\star}))\TailEvent{a^{\star}}+\lambda(a^{\star})\beta T_N(a^{\star},k).
$$
Furthermore,
\begin{align*}
    &N^{\TailExpon{a}}Q_{i,i,N}(k,k)\nonumber\\
    &\quad=N^{\TailExpon{a}}\mathbb{P}\bigg(W_A\big(T_N(a^{\star},k)\big)+\HatBrQueueLength{i}{}{}
    >g_N(a,x,k)+\frac{\sigma_A ^2 \left(\sigma ^2+\sigma_A ^2\right)}{\beta  \sigma ^2}\log N+x\sqrt{\log N}\bigg)\nonumber\\
    &\quad=\int_{-\infty}^{\infty}N^{\TailExpon{a}}\probability{\HatBrQueueLength{i}{}{}>g_N(a,x,k)}\frac{\sqrt{\log N} \exp \left(-\frac{\left(\frac{\sigma_A ^2 \left(\sigma ^2+\sigma_A ^2\right)}{\beta  \sigma ^2}\log N+x\sqrt{\log N}\right)^2}{2 \sigma_A ^2 T_N(a^{\star},k)}\right)}{\sqrt{2 \pi } \sigma_A  \sqrt{T_N(a^{\star},k)}}dx.
\end{align*}
We can simplify this expression further and get that
\begin{align*}
&N^{\TailExpon{a}}\probability{\HatBrQueueLength{i}{}{}>g_N(a,x,k)}\frac{\sqrt{\log N} \exp \left(-\frac{\left(\frac{\sigma_A ^2 \left(\sigma ^2+\sigma_A ^2\right)}{\beta  \sigma ^2}\log N+x\sqrt{\log N}\right)^2}{2 \sigma_A ^2 T_N(a^{\star},k)}\right)}{\sqrt{2 \pi } \sigma_A  \sqrt{T_N(a^{\star},k)}}\nonumber\\
 &\quad=N^{\TailExpon{a}}\exp\left(-\frac{2\beta}{\sigma^2+\sigma_A^2}g_N(a,x,k)\right)
  \frac{\sqrt{\log N} \exp \left(-\frac{\left(\frac{\sigma_A ^2 \left(\sigma ^2+\sigma_A ^2\right)}{\beta  \sigma ^2}\log N+x\sqrt{\log N}\right)^2}{2 \sigma_A ^2 T_N(a^{\star},k)}\right)}{\sqrt{2 \pi } \sigma_A  \sqrt{T_N(a^{\star},k)}}\nonumber\\
   &\quad\LimitN\frac{\beta  \sigma  \exp \left(-\frac{\beta ^2 \sigma ^2 \left(x \left(\sigma ^2+\sigma_A ^2\right)-2 \beta  k \sigma_A
   ^2\right)^2}{\sigma_A ^2 \left(\sigma ^2+\sigma_A ^2\right)^4}\right)}{\sqrt{\pi } \sigma_A  \left(\sigma ^2+\sigma_A ^2\right)}.
\end{align*}
Furthermore, following Lemma \ref{lem: hitting time density}, we have that
\begin{align*}
   f_{\HTime{i}{a^{\star}}{0}}\big(T_N(a^{\star},k)\big)f_{\HTime{j}{a^{\star}}{0}}\big(T_N(a^{\star},l)\big)N^2\log N \LimitN \frac{\beta ^2 \exp\left(-\frac{\beta ^4 k^2}{\left(2 a^{\star} \beta +\sigma ^2\right)^2}\right)}{\sqrt{\pi }
   \left(2 a^{\star} \beta +\sigma ^2\right)}\frac{\beta ^2 \exp\left(-\frac{\beta ^4 l^2}{\left(2 a^{\star} \beta +\sigma ^2\right)^2}\right)}{\sqrt{\pi }
   \left(2 a^{\star} \beta +\sigma ^2\right)}.
\end{align*}
Let $0<\delta<\frac{\beta  \sigma ^4 \epsilon }{2 \sigma_A ^2 \left(\sigma ^2+\sigma_A ^2\right)^2}$ and let 
$$
  h_N(a,x,k)
=\TailEvent{a}-(1-\lambda(a^{\star}))\TailEvent{a^{\star}}+\lambda(a^{\star})\beta T_N(a^{\star},k)-(1+\delta)\left(\frac{\sigma_A ^2 \left(\sigma ^2+\sigma_A ^2\right)}{\beta  \sigma ^2}\log N+x\sqrt{\log N}\right).
$$

From D\k{e}bicki et al.\ \cite[Th. 2.3]{debicki2020exact}, we know that 
\begin{align}\label{eq: debicki limit}
   \probability*{\min(\HatBrQueueLength{i}{}{},\HatBrQueueLength{j}{}{})>x}\exp\left(\frac{2\beta}{\sigma^2/2+\sigma_A^2}x\right)\longrightarrow 0,
\end{align}
as $x\to\infty$. We have that
\begin{align*}
 &N^{\TailExpon{a}}\exp\left(-\frac{2\beta}{\sigma^2/2+\sigma_A^2}h_N(a,x,k)\right)
  \frac{\sqrt{\log N} \exp \left(-\frac{\left(\frac{\sigma_A ^2 \left(\sigma ^2+\sigma_A ^2\right)}{\beta  \sigma ^2}\log N+x\sqrt{\log N}\right)^2}{2 \sigma_A ^2 T_N(a^{\star},k)}\right)}{\sqrt{2 \pi } \sigma_A  \sqrt{T_N(a^{\star},k)}}
   \LimitN 0.
\end{align*}
Thus, when $k\leq l$, then 
\begin{multline*}
    \liminf_{N\to\infty} N^{\TailExpon{a}}f_{\HTime{i}{a^{\star}}{0}}\big(T_N(a^{\star},k)\big)f_{\HTime{j}{a^{\star}}{0}}\big(T_N(a^{\star},l)\big)N^2\log N\left(Q_{i,i,N}(k,k)-\frac{1}{2}Q_{i,j,N}(k,l)\right)^+\nonumber\\
    \geq \frac{\beta ^2 \exp\left(-\frac{\beta ^4 k^2}{\left(2 a^{\star} \beta +\sigma ^2\right)^2}\right)}{\sqrt{\pi }
   \left(2 a^{\star} \beta +\sigma ^2\right)}\frac{\beta ^2 \exp\left(-\frac{\beta ^4 l^2}{\left(2 a^{\star} \beta +\sigma ^2\right)^2}\right)}{\sqrt{\pi }
   \left(2 a^{\star} \beta +\sigma ^2\right)}\frac{\beta  \sigma  \exp \left(-\frac{\beta ^2 \sigma ^2 \left(x \left(\sigma ^2+\sigma_A ^2\right)-2 \beta  k \sigma_A
   ^2\right)^2}{\sigma_A ^2 \left(\sigma ^2+\sigma_A ^2\right)^4}\right)}{\sqrt{\pi } \sigma_A  \left(\sigma ^2+\sigma_A ^2\right)}.
\end{multline*}
The case $k>l$ can be treated analogously. Finally, we have
\begin{align*}
    \int_{-\infty}^{\infty}\int_{-\infty}^{\infty}\int_{-\infty}^{\infty}\frac{\beta ^2 \exp\left(-\frac{\beta ^4 k^2}{\left(2 a^{\star} \beta +\sigma ^2\right)^2}\right)}{\sqrt{\pi }
   \left(2 a^{\star} \beta +\sigma ^2\right)}\frac{\beta ^2 \exp\left(-\frac{\beta ^4 l^2}{\left(2 a^{\star} \beta +\sigma ^2\right)^2}\right)}{\sqrt{\pi }
   \left(2 a^{\star} \beta +\sigma ^2\right)}\frac{\beta  \sigma  \exp \left(-\frac{\beta ^2 \sigma ^2 \left(x \left(\sigma ^2+\sigma_A ^2\right)-2 \beta  k \sigma_A
   ^2\right)^2}{\sigma_A ^2 \left(\sigma ^2+\sigma_A ^2\right)^4}\right)}{\sqrt{\pi } \sigma_A  \left(\sigma ^2+\sigma_A ^2\right)}dxdkdl=1.
\end{align*}
By applying Fatou's lemma, Lemma \ref{lem: liminf a > a star} follows.
\end{proof}
\begin{corollary}
Let $(y_N,N\geq 1)$ be a sequence such that $\liminf_{N\to\infty}y_N/\log N=\infty$, $(W_i,i\leq N)$ be i.i.d.\ Brownian motions with standard deviation $\sigma$, $W_A$ be a Brownian motion with standard deviation $\sigma_A$, for all $i$, $W_i$ and $W_A$ are mutually independent, and $\MaxQueue{}{}$, $\TailExpon{a}$, and $\TailEvent{a}$ are given by Equations \eqref{eq: def max queue length}, \eqref{eq: tail expon}, and \eqref{eq: tailevent}, respectively. Then the tail probability of the steady-state maximum queue length satisfies
$$
\probability{\MaxQueue{}{}{}>y_N}\sim N\probability{\BrQueueLength{i}{}{}>y_N},
$$
as $N\to\infty$.
\end{corollary}\label{cor: yN larger than log N}
\begin{proof}
By using the union bound, we have that $\probability{\MaxQueue{}{}{}>y_N}\leq N\probability{\BrQueueLength{i}{}{}>y_N}$. Furthermore, by using Bonferroni's inequality we obtain that $\probability{\MaxQueue{}{}{}>y_N}\geq N\probability{\BrQueueLength{i}{}{}>y_N}-N^2/2\probability{\BrQueueLength{i}{}{}>y_N,\BrQueueLength{j}{}{}>y_N}$. Now, using the limit in \eqref{eq: debicki limit}, we see that
$$
\limsup_{N\to\infty}\frac{N^2/2\probability{\BrQueueLength{i}{}{}>y_N,\BrQueueLength{j}{}{}>y_N}}{N\probability{\BrQueueLength{i}{}{}>y_N}}\leq \limsup_{N\to\infty}\frac{1}{2}\frac{N\exp\left(-\frac{2\beta}{\sigma^2/2+\sigma_A^2}y_N\right)}{\exp\left(-\frac{2\beta}{\sigma^2+\sigma_A^2}y_N\right)}=0.
$$
The corollary follows.
\end{proof}

\subsection{The case $a=a^{\star}$}\label{subsec: a = a star}
In Section \ref{sec: log asymp}, we showed that we have at least two regimes, namely $0<a<a^{\star}$, and $a\geq a^{\star}$. It turns out, that when we investigate sharper asymptotics, that the case $a=a^{\star}$ deserves special attention. In the present section, we establish that in the case $a=a^{\star}$, $\probability{\MaxQueue{}{}>\TailEvent{a^{\star}}}\sim \frac{1}{2} N^{-\TailExpon{a^{\star}}}$, thus the prefactor is 1/2 instead of 1 as in the case $a>a^{\star}$. 
To make the heuristics given in Section \ref{sec: main results}
rigorous, we proceed by deriving asymptotic lower and upper bounds, in two separate lemmas. As in Section \ref{subsec: a > a star}, we prove that the $\liminf$ converges to the desired limit. We do this in Lemma \ref{lem: liminf a = a star}. The proof of this Lemma is similar to the proof of Lemma \ref{lem: liminf a > a star}. However, the simple union bound $N\probability{\BrQueueLength{i}{}{}>\TailEvent{a^{\star}}}\sim N^{-\TailExpon{a^{\star}}}$ is not tight for $a=a^{\star}$. Thus, we also need to prove that the $\limsup$ is tight. We provide this proof in Lemma \ref{lem: limsup a = a star}.
\begin{lemma}\label{lem: liminf a = a star}
Let $a=a^{\star}$, $(W_i,i\leq N)$ be i.i.d.\ Brownian motions with standard deviation $\sigma$, $W_A$ be a Brownian motion with standard deviation $\sigma_A$, for all $i$, $W_i$ and $W_A$ are mutually independent, and $\MaxQueue{}{}$, $\TailExpon{a}$, and $\TailEvent{a}$ are given by Equations \eqref{eq: def max queue length}, \eqref{eq: tail expon}, and \eqref{eq: tailevent}, respectively, then
\begin{align*}
   \liminf_{N\to\infty} N^{\TailExpon{a^{\star}}}\probability{\MaxQueue{}{}>\TailEvent{a^{\star}}}\geq \frac{1}{2}.
\end{align*}
\end{lemma}
\begin{proof}
First of all, we have the lower bound 
\begin{align*}
    \probability{\MaxQueue{}{}>\TailEvent{a^{\star}}}\geq \probability{\MaxQueue{\HTime{i}{a^{\star}}{r}}{}\mathbbm{1}(\HTime{\wedge}{a^{\star}}{r}<\infty)>\TailEvent{a^{\star}}}.
\end{align*}
As in \eqref{eq: bonf inequality} we can bound this further by Bonferroni's inequality to
\begin{align}
    &N\probability*{\BrQueueLength{i}{\HTime{i}{a^{\star}}{r}}{}\mathbbm{1}(\HTime{i}{a^{\star}}{r}<\infty)>\TailEvent{a^{\star}}}\nonumber\\
    &-\binom{N}{2}\probability*{\min(\BrQueueLength{i}{\HTime{i}{a^{\star}}{r}}{}\mathbbm{1}(\HTime{i}{a^{\star}}{r}<\infty),\BrQueueLength{j}{\HTime{j}{a^{\star}}{r}}{}\mathbbm{1}(\HTime{j}{a^{\star}}{r}<\infty))>\TailEvent{a^{\star}}}\nonumber\\
    &\quad\geq\left(N-\frac{N^2}{2}\probability*{\HTime{j}{a^{\star}}{r}<\infty}\right)\probability*{\BrQueueLength{i}{\HTime{i}{a^{\star}}{r}}{}\mathbbm{1}(\HTime{i}{a^{\star}}{r}<\infty)>\TailEvent{a^{\star}}}\label{subeq: liminf a=a star}.
\end{align}
The last step is true because for independent $X$ and $Y$, $\probability{\min(X,Y\mathbbm{1}(Y<\infty))>z}\leq\probability{Y<\infty}\probability{X>z}$. Since $\probability{\HTime{j}{a^{\star}}{r}<\infty}=\exp(-2(1-\lambda(a^{\star}))\beta r/\sigma^2)/N$, we can simplify the expression in \eqref{subeq: liminf a=a star} to
\begin{align}\label{eq: lower bound liminf a = a star}
    \left(1-\frac{\exp\left(-\frac{2(1-\lambda(a^{\star}))\beta r}{\sigma^2}\right)}{2}\right)N\probability*{\BrQueueLength{i}{\HTime{i}{a^{\star}}{r}}{}\mathbbm{1}(\HTime{i}{a^{\star}}{r}<\infty)>\TailEvent{a^{\star}}}.
\end{align}
Following the proof of Lemma \ref{lem: liminf a > a star} we have that 
\begin{multline*}
    N^{\TailExpon{a^{\star}}}\probability{\HatBrQueueLength{i}{}{}>g_N(a^{\star},x,k)-r}\frac{\sqrt{\log N} \exp \left(-\frac{\left(\frac{\sigma_A ^2 \left(\sigma ^2+\sigma_A ^2\right)}{\beta  \sigma ^2}\log N+x\sqrt{\log N}\right)^2}{2 \sigma_A ^2 T_N(a^{\star},k)}\right)}{\sqrt{2 \pi } \sigma_A  \sqrt{T_N(a^{\star},k)}}\\
   \LimitN\frac{\beta  \sigma  \exp \left(-\frac{\beta ^2 \sigma ^2 \left(x \left(\sigma ^2+\sigma_A ^2\right)-2 \beta  k \sigma_A
   ^2\right)^2}{\sigma_A ^2 \left(\sigma ^2+\sigma_A ^2\right)^4}\right)}{\sqrt{\pi } \sigma_A  \left(\sigma ^2+\sigma_A ^2\right)}\exp\left(\frac{2\beta r}{\sigma^2+\sigma_A^2}\right),
\end{multline*}
when $x<\sigma_A^2\beta k/(\sigma^2+\sigma_A^2)$, and 0 otherwise. Thus, by combining this result with the result from Lemma \ref{lem: hitting time density}, for $x<\sigma_A^2\beta k/(\sigma^2+\sigma_A^2)$,
\begin{multline*}
f_{\HTime{i}{a^{\star}}{r}}\big(T_N(a^{\star},k)\big)N\sqrt{\log N}N^{\TailExpon{a^{\star}}}\probability*{\HatBrQueueLength{i}{}{}>g_N(a^{\star},x,k)-r}\frac{\sqrt{\log N} \exp \left(-\frac{\left(\frac{\sigma_A ^2 \left(\sigma ^2+\sigma_A ^2\right)}{\beta  \sigma ^2}\log N+x\sqrt{\log N}\right)^2}{2 \sigma_A ^2 T_N(a^{\star},k)}\right)}{\sqrt{2 \pi } \sigma_A  \sqrt{T_N(a^{\star},k)}}\\
\LimitN\frac{\beta ^2 \sigma ^2 \exp \left(-\frac{\beta  \left(\beta ^3 k^2 \sigma ^4+2 l \left(\sigma ^2+\sigma_A
   ^2\right)^3\right)}{\left(\sigma ^2+\sigma_A ^2\right)^4}\right)}{\sqrt{\pi } \left(\sigma ^2+\sigma_A^2\right)^2}\frac{\beta  \sigma  \exp \left(-\frac{\beta ^2 \sigma ^2 \left(x \left(\sigma ^2+\sigma_A ^2\right)-2 \beta  k \sigma_A
   ^2\right)^2}{\sigma_A ^2 \left(\sigma ^2+\sigma_A ^2\right)^4}\right)}{\sqrt{\pi } \sigma_A  \left(\sigma ^2+\sigma_A ^2\right)}\exp\left(\frac{2\beta r}{\sigma^2+\sigma_A^2}\right).
\end{multline*}
Observe that the integral
\begin{align*}
    \int_{-\infty}^{\infty}\int_{-\infty}^{\frac{\sigma_A^2\beta k}{\sigma^2+\sigma_A^2}}\frac{\beta ^2 \sigma ^2 \exp \left(-\frac{\beta  \left(\beta ^3 k^2 \sigma ^4+2 r \left(\sigma ^2+\sigma_A
   ^2\right)^3\right)}{\left(\sigma ^2+\sigma_A ^2\right)^4}\right)}{\sqrt{\pi } \left(\sigma ^2+\sigma_A
   ^2\right)^2}\frac{\beta  \sigma  \exp \left(-\frac{\beta ^2 \sigma ^2 \left(x \left(\sigma ^2+\sigma_A ^2\right)-2 \beta  k \sigma_A
   ^2\right)^2}{\sigma_A ^2 \left(\sigma ^2+\sigma_A ^2\right)^4}\right)}{\sqrt{\pi } \sigma_A  \left(\sigma ^2+\sigma_A ^2\right)}\exp\left(\frac{2\beta r}{\sigma^2+\sigma_A^2}\right)dxdk=\frac{1}{2}.
\end{align*}
Now, by applying Fatou's lemma,
\begin{align*}
    \liminf_{N\to\infty}N^{\TailExpon{a^{\star}}}N\probability*{\BrQueueLength{i}{\HTime{i}{a^{\star}}{r}}{}\mathbbm{1}(\HTime{i}{a^{\star}}{r}<\infty)>\TailEvent{a^{\star}}}\geq \frac{1}{2},
\end{align*}
and thus by applying this on the expression in \eqref{eq: lower bound liminf a = a star}, we get that
\begin{align*}
       \liminf_{N\to\infty} N^{\TailExpon{a^{\star}}}\probability{\MaxQueue{}{}>\TailEvent{a^{\star}}}\geq \frac{1}{2}\left(1-\frac{\exp\left(-\frac{2(1-\lambda(a^{\star}))\beta r}{\sigma^2}\right)}{2}\right)\overset{r\to\infty}{\longrightarrow} \frac{1}{2}.
\end{align*}
\end{proof}
\begin{lemma}\label{lem: limsup a = a star}
Let $a=a^{\star}$, $(W_i,i\leq N)$ be i.i.d.\ Brownian motions with standard deviation $\sigma$, $W_A$ be a Brownian motion with standard deviation $\sigma_A$, for all $i$, $W_i$ and $W_A$ are mutually independent, and $\MaxQueue{}{}$, $\TailExpon{a}$, and $\TailEvent{a}$ are given by Equations \eqref{eq: def max queue length}, \eqref{eq: tail expon}, and \eqref{eq: tailevent}, respectively, then
\begin{align*}
   \limsup_{N\to\infty} N^{\TailExpon{a^{\star}}}\probability{\MaxQueue{}{}>\TailEvent{a^{\star}}}\leq \frac{1}{2}.
\end{align*}
\end{lemma}
\begin{proof}
Let $\HTimeA{a^{\star}}{r}=\inf\{t:W_A(t)-\lambda(a^{\star})\beta t>\lambda(a^{\star})\TailEvent{a^{\star}}+r\}$. Following Equation \eqref{eq: repres max Q hitting times} and the upper bound in \eqref{eq: upper lower bounds repr max Q hitting times}, we have that
\begin{align}\label{subeq: upper bound a = a star}
   \probability{\MaxQueue{}{}>\TailEvent{a^{\star}}}
   \leq\probability{\MaxQueue{\HTimeA{a^{\star}}{r}}{}\mathbbm{1}(\HTimeA{a^{\star}}{r}<\infty)>\TailEvent{a^{\star}}}
   +\probability{\MaxQueue{\HTime{i}{a^{\star}}{-r}}{}\mathbbm{1}(\HTime{\wedge}{a^{\star}}{-r}<\infty)>\TailEvent{a^{\star}}}.
\end{align}
Observe that we can bound the first term in \eqref{subeq: upper bound a = a star} as
\begin{align}\label{eq: a = a star upper bound 1}
    \probability{\MaxQueue{\HTimeA{a^{\star}}{r}}{}\mathbbm{1}(\HTimeA{a^{\star}}{r}<\infty)>\TailEvent{a^{\star}}}\leq   \probability{\HTimeA{a^{\star}}{r}<\infty}=N^{-\TailExpon{a^{\star}}}\exp\left(-\frac{2\lambda(a^{\star})\beta r}{\sigma_A^2}\right),
\end{align}
and the second term in \eqref{subeq: upper bound a = a star} as
\begin{align}
   &N^{\TailExpon{a^{\star}}}\probability{\MaxQueue{\HTime{i}{a^{\star}}{-r}}{}\mathbbm{1}(\HTime{\wedge}{a^{\star}}{-r}<\infty)>\TailEvent{a^{\star}}}\nonumber\\
   &\quad\leq N^{\TailExpon{a^{\star}}}N\probability*{\BrQueueLength{i}{\HTime{i}{a^{\star}}{-r}}{}\mathbbm{1}(\HTime{i}{a^{\star}}{-r}<\infty)>\TailEvent{a^{\star}}}\nonumber\\
   &\quad=\int_{-\infty}^{\infty}N^{\TailExpon{a^{\star}}}N\probability*{\BrQueueLength{i}{\HTime{i}{a^{\star}}{-r}}{}>\TailEvent{a^{\star}}\Big|\HTime{i}{a^{\star}}{-r}=T_N(a^{\star},k)}
   f_{\HTime{i}{a^{\star}}{-r}}\big(T_N(a^{\star},k)\big)\sqrt{\log N}dk.\label{eq: integrand a= a star limsup}
\end{align}
Now, we examine the parts of the integrand of this integral separately. First, note that, following Definition \ref{def: hitting time},
\begin{multline*}
    \probability*{\BrQueueLength{i}{\HTime{i}{a^{\star}}{-r}}{}>\TailEvent{a^{\star}}\Big|\HTime{i}{a^{\star}}{-r}=T_N(a^{\star},k)}\\
    =\probability*{W_A(\HTime{i}{a^{\star}}{-r})+\HatBrQueueLength{i}{}{}>\lambda(a^{\star})\TailEvent{a^{\star}}+r+\lambda(a^{\star})\beta \HTime{i}{a^{\star}}{-r}\Big|\HTime{i}{a^{\star}}{-r}=T_N(a^{\star},k)},
\end{multline*}
We can analyze this probability using Lemma \ref{lem: convolution normal exponential} by taking $x_N=2\lambda(a^{\star})\TailEvent{a^{\star}}+\lambda(a^{\star})\beta k\sqrt{\log N}+r$, $\eta_N=\sigma_A\sqrt{T_N(a^{\star},k)}$, and $\mu=2\beta/(\sigma^2+\sigma_A^2)$. Write 
\begin{align*}
    \frac{x_N-\mu\eta_N^2}{\sqrt{2}\eta_N}=&\frac{2\lambda(a^{\star})\TailEvent{a^{\star}}+\lambda(a^{\star})\beta k\sqrt{\log N}+r-\frac{2\beta}{\sigma^2+\sigma_A^2}\sigma_A^2T_N(a^{\star},k)}{\sqrt{2}\sqrt{\sigma_A^2T_N(a^{\star},k)}}\nonumber\\
    =&\frac{r-\lambda(a^{\star})\beta k\sqrt{\log N}}{\sqrt{2}\sqrt{\sigma_A^2T_N(a^{\star},k)}}\LimitN -\frac{\sqrt{2} \beta ^3 k \sigma_A  \sigma ^2}{\left(\sigma ^2+\sigma_A ^2\right)^3}.
\end{align*}
The first term in \eqref{eq: convolution tail 1} of Lemma \ref{lem: convolution normal exponential} satisfies
\begin{align*}
    \frac{\sigma  \exp \left(-\frac{\beta  \left(\beta ^3 k^2 \sigma_A ^2 \sigma ^2+2 r \left(\sigma
   ^2+\sigma_A ^2\right)^3\right)}{\left(\sigma ^2+\sigma_A ^2\right)^4}\right)}{2 \sqrt{\pi } \sigma_A
   }\frac{N^{-\TailExpon{a^{\star}}}}{\sqrt{\log N}},
\end{align*}
and the second term satisfies
\begin{align*}
   \frac{1}{2} e^{\frac{1}{2} \mu  \left(\mu  \eta_N ^2-2 x_N\right)}\bigg(1+\text{erf}\bigg(-\frac{\sqrt{2} \beta ^3 k \sigma_A  \sigma ^2}{\left(\sigma ^2+\sigma_A ^2\right)^3}\bigg)\bigg)\sim \frac{1}{2}\exp\left(-\frac{2 \beta  r}{\sigma ^2+\sigma_A ^2}\right)\bigg(1+\text{erf}\bigg(-\frac{\sqrt{2} \beta ^3 k \sigma_A  \sigma ^2}{\left(\sigma ^2+\sigma_A ^2\right)^3}\bigg)\bigg)N^{-\TailExpon{a^{\star}}},
\end{align*}
as $N\to\infty$. So, we can conclude that 
$$
\probability*{\BrQueueLength{i}{\HTime{i}{a^{\star}}{-r}}{}>\TailEvent{a^{\star}}\Big| \HTime{i}{a^{\star}}{-r}=T_N(a^{\star},k)}
  \sim \frac{1}{2}\exp\left(-\frac{2 \beta  r}{\sigma ^2+\sigma_A ^2}\right)\bigg(1+\text{erf}\bigg(-\frac{\sqrt{2} \beta ^3 k \sigma_A  \sigma ^2}{\left(\sigma ^2+\sigma_A ^2\right)^3}\bigg)\bigg)N^{-\TailExpon{a^{\star}}},
$$
as $N\to\infty$. Second, following Lemma \ref{lem: hitting time density}, the density of the hitting time $\HTime{i}{a^{\star}}{-r}$ appears in the integrand in \eqref{eq: integrand a= a star limsup}, and satisfies
\begin{align*}
    Nf_{\HTime{i}{a^{\star}}{-r}}\big(T_N(a^{\star},k)\big)\sqrt{\log N}\LimitN& \frac{\beta ^2 \exp \left(\frac{\beta  \left(8 a^{{\star}^2} \beta ^2 r-\beta ^3 k^2 \sigma  \sqrt{2 a^{\star} \beta +\sigma ^2}+8 a^{\star} \beta 
   r \sigma ^2+2 r \sigma ^4\right)}{\sigma  \left(2 a^{\star} \beta +\sigma ^2\right)^{5/2}}\right)}{\sqrt{\pi } \left(2 a^{\star} \beta
   +\sigma ^2\right)}\nonumber\\
   =& \frac{\beta ^2 \sigma ^2 \exp \left(\frac{\beta  \left(2 r \left(\sigma ^2+\sigma_A
   ^2\right)^3-\beta ^3 k^2 \sigma ^4\right)}{\left(\sigma ^2+\sigma_A
   ^2\right)^4}\right)}{\sqrt{\pi } \left(\sigma ^2+\sigma_A ^2\right)^2}.
\end{align*}
Thus, for the integrand in \eqref{eq: integrand a= a star limsup} we have that
\begin{multline*}
    N^{\TailExpon{a^{\star}}}N\probability*{\BrQueueLength{i}{\HTime{i}{a^{\star}}{-r}}{}>\TailEvent{a^{\star}}\Big|\HTime{i}{a^{\star}}{-r}=T_N(a^{\star},k)}
    f_{\HTime{i}{a^{\star}}{-r}}\big(T_N(a^{\star},k)\big)\sqrt{\log N}\\
    \LimitN \frac{\beta ^2 \sigma ^2 \bigg(1+\text{erf}\bigg(-\frac{\sqrt{2} \beta ^3 k \sigma_A  \sigma ^2}{\left(\sigma ^2+\sigma_A ^2\right)^3}\bigg)\bigg) \exp \left(\frac{\beta  \left(2 r \left(\sigma ^2+\sigma_A ^2\right)^3-\beta ^3 k^2 \sigma
   ^4\right)}{\left(\sigma ^2+\sigma_A ^2\right)^4}-\frac{2 \beta  r}{\sigma ^2+\sigma_A ^2}\right)}{2 \sqrt{\pi } \left(\sigma
   ^2+\sigma_A ^2\right)^2}.
\end{multline*}
When we integrate this result we get
\begin{align*}
    \int_{-\infty}^{\infty}\frac{\beta ^2 \sigma ^2 \bigg(1+\text{erf}\bigg(-\frac{\sqrt{2} \beta ^3 k \sigma_A  \sigma ^2}{\left(\sigma ^2+\sigma_A ^2\right)^3}\bigg)\bigg) \exp \left(\frac{\beta  \left(2 r \left(\sigma ^2+\sigma_A ^2\right)^3-\beta ^3 k^2 \sigma
   ^4\right)}{\left(\sigma ^2+\sigma_A ^2\right)^4}-\frac{2 \beta  r}{\sigma ^2+\sigma_A ^2}\right)}{2 \sqrt{\pi } \left(\sigma
   ^2+\sigma_A ^2\right)^2}dk=\frac{1}{2}.
\end{align*}
Now, because
\begin{align*}
    &N^{\TailExpon{a^{\star}}}\probability*{\BrQueueLength{i}{\HTime{i}{a^{\star}}{-r}}{}>\TailEvent{a^{\star}}\Big|\HTime{i}{a^{\star}}{-r}=T_N(a^{\star},k)}\nonumber\\
    &\quad\leq N^{\TailExpon{a^{\star}}}\probability*{\sup_{s>0}\left(W_A(s)-\lambda(a^{\star})\beta s\right)>\lambda(a^{\star})\TailEvent{a^{\star}}+r}\nonumber\\
    &\quad = N^{\TailExpon{a^{\star}}}\exp\left(-\frac{2\lambda(a^{\star})\beta}{\sigma_A^2}(\lambda(a^{\star})\TailEvent{a^{\star}}+r)\right)=\exp\left(-\frac{2\lambda(a^{\star})\beta r}{\sigma_A^2}\right),
\end{align*}
and 
\begin{align*}
    \lim_{N\to\infty}\int_{-\infty}^{\infty}Nf_{\HTime{i}{a^{\star}}{-r}}\big(T_N(a^{\star},k)\big)\sqrt{\log N}dk= \int_{-\infty}^{\infty}\lim_{N\to\infty}Nf_{\HTime{i}{a^{\star}}{-r}}\big(T_N(a^{\star},k)\big)\sqrt{\log N}dk,
\end{align*}
we can use Lemma \ref{lem: conv sequence integrals} to conclude that
\begin{align}\label{eq: a = a star upper bound 2}
   \limsup_{N\to\infty} N^{\TailExpon{a^{\star}}}N\probability*{\BrQueueLength{i}{\HTime{i}{a^{\star}}{-r}}{}\mathbbm{1}(\HTime{i}{a^{\star}}{-r}<\infty)>\TailEvent{a^{\star}}}\leq\frac{1}{2}.
\end{align}
Now, after combining the bounds in \eqref{eq: a = a star upper bound 1} and \eqref{eq: a = a star upper bound 2}, 
\begin{align*}
   \limsup_{N\to\infty} N^{\TailExpon{a^{\star}}}\probability{\MaxQueue{}{}>\TailEvent{a^{\star}}}\leq \frac{1}{2}+\exp\left(-\frac{2\lambda(a^{\star})\beta r}{\sigma_A^2}\right)\overset{r\to\infty}{\longrightarrow}\frac{1}{2}.
\end{align*}
\end{proof}

\subsection{The case $0<a<a^{\star}$}\label{subsec: 0<a< a star}
As we have proven the exact asymptotics for the cases $a>a^{\star}$ and $a=a^{\star}$ in Theorems \ref{thm: exact asymp a > a star} and \ref{thm: exact asymp a = a star}, respectively, we now turn to the proof of Theorem \ref{thm: exact asymp 0<a< a star}. In Theorem \ref{thm: log asymp} we have shown that $\TailExpon{a}=\frac{2a\beta+2\sigma^2-2\sigma\sqrt{2a\beta+\sigma^2}}{\sigma_A^2}$, thus we expect highly dependent behavior because this indicates that the union upper bound $\probability{\MaxQueue{}{}>\TailEvent{a}}\leq N\probability{\BrQueueLength{i}{}{}>\TailEvent{a}}$ is not sharp when $0<a<a^{\star}$, as is explained in the proof of Lemma \ref{lem: limsup log asymp}.
\begin{proof}[Proof of Theorem \ref{thm: exact asymp 0<a< a star}]
First of all, we prove Equation \eqref{eq: 0<a< a star liminf}. We write $r_N=\frac{\sigma  \sqrt{2 a \beta +\sigma ^2}}{4 \beta }\log\log N$. Let $\HTimeA{a}{r_N}=\inf\{t>0:W_A(t)-\lambda(a)\beta t>\lambda(a)\TailEvent{a}+r_N\}$. Let $f_{\HTimeA{a}{r_N}}$ be its density.
Observe that
\begin{align}\label{eq: 0<a<a star integrand}
    &\probability{\MaxQueue{}{}>\TailEvent{a}}\nonumber\\
    &\quad\geq  \probability{\MaxQueue{\HTimeA{a}{r_N}}{\HTimeA{a}{r_N}}\mathbbm{1}(\HTimeA{a}{r_N}<\infty)>\TailEvent{a}}\nonumber\\
    &\quad= \int_{-\infty}^{\infty}\probability*{\max_{i\leq N}W_i\big(T_N(a,k)\big)-(1-\lambda(a))\beta T_N(a,k)>(1-\lambda(a))\TailEvent{a}-r_N}
     f_{\HTimeA{a}{r_N}}\big(T_N(a,k)\big)\sqrt{\log N}dk.
\end{align}
As in the proof of Lemma \ref{lem: limsup a = a star}, we analyze the components of the integrand of \eqref{eq: 0<a<a star integrand} separately. Following a similar derivation as in Lemma \ref{lem: hitting time density}, we see that the term $f_{\HTimeA{a}{r_N}}\big(T_N(a,k)\big)\sqrt{\log N}$ in \eqref{eq: 0<a<a star integrand} satisfies
\begin{multline}\label{eq: hitting time density 0<a<a star convergence}
    N^{\TailExpon{a}}(\log N)^{\frac{\lambda(a)}{1-\lambda(a)}\frac{\sigma^2}{2\sigma_A^2}} f_{\HTimeA{a}{r_N}}\big(T_N(a,k)\big)\sqrt{\log N}\\
    \LimitN \frac{\beta ^2 \left(\sigma  \left(\sigma -\sqrt{2 a \beta +\sigma ^2}\right)+2 a \beta \right) \exp \left(-\frac{\beta
   ^4 k^2 \left(\sigma -\sqrt{2 a \beta +\sigma ^2}\right)^2}{\sigma_A ^2 \left(2 a \beta +\sigma
   ^2\right)^2}\right)}{\sqrt{\pi } \sigma_A  \left(2 a \beta +\sigma ^2\right)^{3/2}}.
\end{multline}
Moreover, a result in extreme value theory states that when  $b_N=\sqrt{2\log N}-\log(4\pi\log N)/(2\sqrt{2\log N})$, then
\begin{align*}
b_N\left(\frac{\max_{i\leq N}W_i(d\log N)}{\sigma\sqrt{d\log N}}-b_N\right)\LimitD G,
\end{align*}
with $G\overset{d}{=}\text{Gumbel}$, as $N\to\infty$, cf.\ \cite[p.\ 11, Ex.\ 1.1.7]{de2007extreme} for a proof. From this it follows that the term $\probability*{\max_{i\leq N}W_i\big(T_N(a,k)\big)-(1-\lambda(a))\beta T_N(a,k)>(1-\lambda(a))\TailEvent{a}-r_N}$ in \eqref{eq: 0<a<a star integrand} satisfies
\begin{multline}\label{eq: convergence gumbel}
\probability*{\max_{i\leq N}W_i\big(T_N(a,k)\big)-(1-\lambda(a))\beta T_N(a,k)>(1-\lambda(a))\TailEvent{a}-r_N}\\
    \LimitN 1-\exp\left(-\frac{\exp\bigg(-\frac{\beta ^4 k^2}{\left(2 a \beta +\sigma ^2\right)^2}\bigg)}{2 \sqrt{\pi }}\right).    
\end{multline}
Thus, the product of the limits in \eqref{eq: hitting time density 0<a<a star convergence} and \eqref{eq: convergence gumbel} gives the tail asymptotics of the integrand in \eqref{eq: 0<a<a star integrand}. Now, by applying Fatou's lemma, we obtain a sharper than logarithmic lower bound on the asymptotics for the maximum queue length, and is given in \eqref{eq: 0<a< a star liminf}. \\

In order to prove \eqref{eq: 0<a< a star limsup}, we use the upper bound given in \eqref{eq: upper lower bounds repr max Q hitting times} and observe that 
\begin{align}
    \probability{\MaxQueue{}{}>\TailEvent{a}}
    \leq&\probability{\MaxQueue{\HTimeA{a}{r_N}}{}\mathbbm{1}(\HTimeA{a}{r_N}<\infty)>\TailEvent{a}}\label{subeq: limsup 0<a<a star 1}\\
   &\quad+\probability{\MaxQueue{\HTime{i}{a}{-r_N}}{}\mathbbm{1}(\HTime{\wedge}{a}{-r_N}<\infty)>\TailEvent{a}}\label{subeq: limsup 0<a<a star 2}.
\end{align}
We can bound the expression in \eqref{subeq: limsup 0<a<a star 1} as follows:
\begin{align}
   \probability{\MaxQueue{\HTimeA{a}{r_N}}{}\mathbbm{1}(\HTimeA{a}{r_N}<\infty)>\TailEvent{a}}\leq \probability{\HTimeA{a}{r_N}<\infty} =N^{-\TailExpon{a}}(\log N)^{-\frac{\lambda(a)}{1-\lambda(a)}\frac{\sigma^2}{2\sigma_A^2}}.
\end{align}
Therefore,
$$
\limsup_{N\to\infty}N^{\TailExpon{a}}(\log N)^{\frac{\lambda(a)}{1-\lambda(a)}\frac{\sigma^2}{2\sigma_A^2}} \probability{\MaxQueue{\HTimeA{a}{r_N}}{}\mathbbm{1}(\HTimeA{a}{r_N}<\infty)>\TailEvent{a}}\leq 1.
$$
Thus, because of the bounds given in \eqref{subeq: limsup 0<a<a star 1} and \eqref{subeq: limsup 0<a<a star 2}, to prove that \eqref{eq: 0<a< a star limsup} holds, it is left to show that
$$\limsup_{N\to\infty}N^{\TailExpon{a}}(\log N)^{\frac{\lambda(a)}{1-\lambda(a)}\frac{\sigma^2}{2\sigma_A^2}}\probability{\MaxQueue{\HTime{i}{a}{-r_N}}{}\mathbbm{1}(\HTime{\wedge}{a}{-r_N}<\infty)>\TailEvent{a}}<\infty.$$ 
To prove this, observe that, by using the union bound and by conditioning on the hitting time $\HTime{i}{a}{-r_N}$ the expression in \eqref{subeq: limsup 0<a<a star 2} satisfies
\begin{align}
    &\probability{\MaxQueue{\HTime{i}{a}{-r_N}}{}\mathbbm{1}(\HTime{\wedge}{a}{-r_N}<\infty)>\TailEvent{a}}\nonumber\\
    &\quad\leq N\probability{\BrQueueLength{i}{\HTime{i}{a}{-r_N}}{}\mathbbm{1}(\HTime{i}{a}{-r_N}<\infty)>\TailEvent{a}}\nonumber\\
    &\quad=\int_{-\infty}^{\infty}N\probability{\BrQueueLength{i}{\HTime{i}{a}{-r_N}}{}>\TailEvent{a}\mid\HTime{i}{a}{-r_N}=T_N(a,k)}f_{\HTime{i}{a}{-r_N}}\big(T_N(a,k)\big)\sqrt{\log N}dk\label{subeq: integral 0<a<a star upper bound}.
\end{align}
Now, we can use Lemma \ref{lem: conv sequence integrals} to show convergence of the integral in \eqref{subeq: integral 0<a<a star upper bound}. Following a similar analysis as in Lemma \ref{lem: hitting time density}, we have that 
\begin{align*}
    N\frac{1}{\sqrt{\log N}}\sqrt{\log N}f_{\HTime{i}{a}{-r_N}}\big(T_N(a,k)\big)\LimitN \frac{\beta ^2 \exp\left(-\frac{\beta ^4 k^2}{\left(2 a \beta +\sigma ^2\right)^2}\right)}{\sqrt{\pi } \left(2 a \beta +\sigma
   ^2\right)}.
\end{align*}
Furthermore,
\begin{align*}
    \int_{-\infty}^{\infty}\frac{\beta ^2 e^{-\frac{\beta ^4 k^2}{\left(2 a \beta +\sigma ^2\right)^2}}}{\sqrt{\pi } \left(2 a \beta +\sigma
   ^2\right)}dk=\int_{-\infty}^{\infty}\frac{N}{\sqrt{\log N}}\sqrt{\log N}f_{\HTime{i}{a}{-r_N}}\big(T_N(a,k)\big)dk=1.
\end{align*}
Thus, the first and second condition in Lemma \ref{lem: conv sequence integrals} hold. 
Thus, we now only need to analyze 
\begin{multline}\label{eq: expon plus normal 0<a< a star}
    \probability{\BrQueueLength{i}{\HTime{i}{a}{-r_N}}{}>\TailEvent{a}\mid\HTime{i}{a}{-r_N}=T_N(a,k)}\\
    =\probability{W_A(\HTime{i}{a}{-r_N})+\HatBrQueueLength{i}{}{}>\lambda(a)\TailEvent{a}+r_N+\lambda(a)\beta \HTime{i}{a}{-r}\mid\HTime{i}{a}{-r_N}=T_N(a,k)},
\end{multline}
which is a component in the integrand in \eqref{subeq: integral 0<a<a star upper bound}. We show that this expression satisfies the third and fourth condition of Lemma \ref{lem: conv sequence integrals}, by proving pointwise convergence and by proving that this probability is uniformly bounded by a constant. To do this, first observe that the random variable in \eqref{eq: expon plus normal 0<a< a star} has the form of the sum of a normally distributed random variable and an exponentially distributed random variable, hence we can follow the framework of Lemma \ref{lem: convolution normal exponential} in order to analyze this probability, we take $x_N=2\lambda(a)\TailEvent{a}+\lambda(a)\beta k\sqrt{\log N}+r_N$, $\eta_N=\sigma_A\sqrt{T_N(a,k)}$, and $\mu=2\beta/(\sigma^2+\sigma_A^2)$. Now, the expression in \eqref{eq: expon plus normal 0<a< a star} can be written in the form of Equation \eqref{eq: exact expression convolution normal exponential}. Furthermore, observe that
\begin{align*}
    \frac{x_N-\mu\eta_N^2}{\sqrt{2}\eta_N}=&\frac{2\lambda(a)\TailEvent{a}+\lambda(a)\beta k\sqrt{\log N}+r_N-\frac{2\beta}{\sigma^2+\sigma_A^2}\sigma_A^2T_N(a,k)}{\sqrt{2}\sqrt{\sigma_A^2T_N(a,k)}}\LimitN -\infty.
\end{align*}
Thus, for $0<a<a^{\star}$, we are in the third situation of Lemma \ref{lem: convolution normal exponential}. The first term in \eqref{eq: convolution tail 3} satisfies
$$
    \frac{\eta_N  e^{-\frac{x_N^2}{2 \eta_N ^2}}}{\sqrt{2 \pi } x_N}\sim\frac{\sigma_A  \exp \left(-\frac{\beta ^4 k^2 \left(\sigma -\sqrt{2 a \beta +\sigma ^2}\right)^2}{\sigma_A ^2 \left(2 a \beta +\sigma ^2\right)^2}\right)}{2 \sqrt{\pi } \left(\sqrt{2 a \beta +\sigma ^2}-\sigma
   \right)}\log N^{-{\frac{\lambda(a)}{1-\lambda(a)}\frac{\sigma^2}{2\sigma_A^2}}}N^{-\TailExpon{a}}\frac{1}{\sqrt{\log N}},
$$
as $N\to\infty$. Furthermore, we have for all $t>0$ that 
$$
    \probability*{W_A(t)-\lambda(a)\beta t>x}\leq \probability*{W_A(x/(\lambda(a)\beta))>2x}.
$$
From this it follows that the first part in \eqref{eq: exact expression convolution normal exponential} satisfies
\begin{align*}
&\probability*{\eta_NX>x_N}\\
    &\quad=\probability*{W_A(\HTime{i}{a}{-r_N})>\lambda(a)\TailEvent{a}+r_N+\lambda(a)\beta \HTime{i}{a}{-r}\bigg|\HTime{i}{a}{-r_N}=T_N(a,k)}\nonumber\\
    &  \quad\leq \probability*{W_A(\HTime{i}{a}{-r_N})>\lambda(a)\TailEvent{a}+r_N+\lambda(a)\beta \HTime{i}{a}{-r}\bigg|\HTime{i}{a}{-r_N}=\frac{\TailEvent{a}}{\beta}+\frac{r_N}{\lambda(a)\beta}}\nonumber\\
    &\quad\sim\frac{\sigma_A }{2 \sqrt{\pi } \left(\sqrt{2 a \beta +\sigma ^2}-\sigma
   \right)}\log N^{-{\frac{\lambda(a)}{1-\lambda(a)}\frac{\sigma^2}{2\sigma_A^2}}}N^{-\TailExpon{a}}\frac{1}{\sqrt{\log N}},
\end{align*}
as $N\to\infty$. So there exists an $\epsilon>0$ and an $N_{\epsilon}$ such that for $N>N_{\epsilon}$ and all $k>-\TailEvent{a}/(\beta\sqrt{\log N})$,
\begin{multline}\label{eq: upper bound 0<a< a star convergence integral}
   (\log N)^{{\frac{\lambda(a)}{1-\lambda(a)}\frac{\sigma^2}{2\sigma_A^2}}}N^{\TailExpon{a}}\sqrt{\log N} \probability*{W_A(\HTime{i}{a}{-r_N})>\lambda(a)\TailEvent{a}+r_N+\lambda(a)\beta \HTime{i}{a}{-r_N}\bigg|\HTime{i}{a}{-r_N}=T_N(a,k)}\\
   \leq \frac{\sigma_A }{2 \sqrt{\pi } \left(\sqrt{2 a \beta +\sigma ^2}-\sigma
   \right)}+\epsilon.
\end{multline}
The second term in \eqref{eq: convolution tail 3} satisfies
\begin{multline}\label{eq: limsup 0<a< a star second part}
-\frac{1}{\sqrt{2\pi}}e^{\frac{1}{2} \mu  \left(\mu  \eta_N ^2-2 x_N\right)}\frac{ \eta_N  e^{-\frac{\left(x_N-\mu  \eta_N ^2\right)^2}{2 \eta_N ^2}}}{x_N-\mu 
   \eta_N ^2}\\
   \sim \frac{\sigma_A  \left(\sigma ^2+\sigma_A ^2\right) \exp \left(-\frac{2 \beta ^4 k^2 \left(\sigma ^2 \left(\sqrt{2 a \beta +\sigma ^2}-\sigma \right)+a \beta  \left(\sqrt{2 a \beta +\sigma ^2}-2 \sigma
   \right)\right)}{\sigma_A ^2 \left(2 a \beta +\sigma ^2\right)^{5/2}}\right)}{2 \sqrt{\pi } \sigma  \left(\sigma  \left(\sigma -\sqrt{2 a \beta +\sigma ^2}\right)+\sigma_A ^2\right)}\log N^{-{\frac{\lambda(a)}{1-\lambda(a)}\frac{\sigma^2}{2\sigma_A^2}}}N^{-\TailExpon{a}}\frac{1}{\sqrt{\log N}},
\end{multline}
as $N\to\infty$. In this case first observe that in Equation \eqref{eq: exact expression convolution normal exponential} the exact expression of the convolution term equals
\begin{align*}
    \int_{-\infty}^{x_N/\eta_N}\probability*{\frac{1}{\mu}E>x_N-\eta_Nz}\frac{e^{-\frac{z^2}{2}}}{\sqrt{2 \pi }}dz=\frac{1}{2}
   \left(\text{erf}\left(\frac{x_N-\mu  \eta_N ^2}{\sqrt{2} \eta_N }\right)+1\right)e^{\frac{1}{2} \mu  \left(\mu  \eta_N ^2-2 x_N\right)}.
\end{align*}
Second, observe that this can be further rewritten into
\begin{align*}
    &\frac{1}{2}
   \left(\text{erf}\left(\frac{x_N-\mu  \eta_N ^2}{\sqrt{2} \eta_N }\right)+1\right)e^{\frac{1}{2} \mu  \left(\mu  \eta_N ^2-2 x_N\right)}\nonumber\\
   &\quad=\probability*{W_A(\HTime{i}{a}{-r_N})>\frac{2\beta}{\sigma^2+\sigma_A^2}\sigma_A^2\HTime{i}{a}{-r_N}-\lambda(a)\TailEvent{a}-r_N-\lambda(a)\beta \HTime{i}{a}{-r_N}\Big|\HTime{i}{a}{-r_N}=T_N(a,k)}\nonumber\\
   &\quad\quad\cdot\exp\left(\frac{1}{2}\frac{2\beta}{\sigma^2+\sigma_A^2}\left(\frac{2\beta}{\sigma^2+\sigma_A^2}\sigma_A^2T_N(a,k)-2\lambda(a)\TailEvent{a}-2\lambda(a)\beta T_N(a,k)-2r_N\right)\right).
\end{align*}
Thus, the expression that we are investigating is a product of a tail probability of a Gaussian random variable and an exponential function. With an analogous derivation as for the first term in \eqref{eq: convolution tail 3}, due to the expression in \eqref{eq: limsup 0<a< a star second part} we can bound for all $t>0$
\begin{multline*}
    (\log N)^{{\frac{\lambda(a)}{1-\lambda(a)}\frac{\sigma^2}{2\sigma_A^2}}}N^{\TailExpon{a}}\sqrt{\log N}\probability*{W_A(t)>\frac{2\beta}{\sigma^2+\sigma_A^2}\sigma_A^2t-\lambda(a)\TailEvent{a}-r_N-\lambda(a)\beta t}\\
    \cdot\exp\left(\frac{1}{2}\frac{2\beta}{\sigma^2+\sigma_A^2}\left(\frac{2\beta}{\sigma^2+\sigma_A^2}\sigma_A^2t-2\lambda(a)\TailEvent{a}-2\lambda(a)\beta t-2r_N\right)\right).
    \end{multline*}
Hence, due to this and the upper bound given in \eqref{eq: upper bound 0<a< a star convergence integral}, we have that the third and fourth condition of Lemma \ref{lem: conv sequence integrals} are satisfied. Thus, in the end we know that
\begin{align*}
    &(\log N)^{{\frac{\lambda(a)}{1-\lambda(a)}\frac{\sigma^2}{2\sigma_A^2}}}N^{\TailExpon{a}}N\probability*{\BrQueueLength{i}{\HTime{i}{a}{-r_N}}{}>\TailEvent{a}\Big|\HTime{i}{a}{-r_N}=T_N(a,k)}f_{\HTime{i}{a}{-r_N}}\big(T_N(a,k)\big)\sqrt{\log N}\nonumber\\
    \LimitN& \left(\frac{\sigma_A  \exp \left(-\frac{\beta ^4 k^2 \left(\sigma -\sqrt{2 a \beta +\sigma ^2}\right)^2}{\sigma_A ^2 \left(2 a \beta +\sigma ^2\right)^2}\right)}{2 \sqrt{\pi } \left(\sqrt{2 a \beta +\sigma ^2}-\sigma
   \right)}+\frac{\sigma_A  \left(\sigma ^2+\sigma_A ^2\right) \exp \left(-\frac{2 \beta ^4 k^2 \left(\sigma ^2 \left(\sqrt{2 a \beta +\sigma ^2}-\sigma \right)+a \beta  \left(\sqrt{2 a \beta +\sigma ^2}-2 \sigma
   \right)\right)}{\sigma_A ^2 \left(2 a \beta +\sigma ^2\right)^{5/2}}\right)}{2 \sqrt{\pi } \sigma  \left(\sigma  \left(\sigma -\sqrt{2 a \beta +\sigma ^2}\right)+\sigma_A ^2\right)}\right)\nonumber\\
    &\cdot \frac{\beta ^2 e^{-\frac{\beta ^4 k^2}{\left(2 a \beta +\sigma ^2\right)^2}}}{\sqrt{\pi } \left(2 a \beta +\sigma
   ^2\right)},
\end{align*}
and we apply Lemma \ref{lem: conv sequence integrals} to conclude that \eqref{eq: 0<a< a star limsup} holds.
\end{proof}
\begin{remark}
We have stated in Theorem \ref{thm: exact asymp 0<a< a star} that we can prove a lower and upper bound which are sharper than logarithmic, however we do not specify these bounds, but from the proof of Theorem \ref{thm: exact asymp 0<a< a star} it becomes clear that 
\begin{align*}
    &\liminf_{N\to\infty}N^{\TailExpon{a}}(\log N)^{\frac{\lambda(a)}{1-\lambda(a)}\frac{\sigma^2}{2\sigma_A^2}}\probability{\MaxQueue{}{}>\TailEvent{a}}\nonumber\\
    \geq&\int_{-\infty}^{\infty}\frac{\beta ^2 \left(\sigma  \left(\sigma -\sqrt{2 a \beta +\sigma ^2}\right)+2 a \beta \right) \exp \left(-\frac{\beta
   ^4 k^2 \left(\sigma -\sqrt{2 a \beta +\sigma ^2}\right)^2}{\sigma_A ^2 \left(2 a \beta +\sigma
   ^2\right)^2}\right)}{\sqrt{\pi } \sigma_A  \left(2 a \beta +\sigma ^2\right)^{3/2}} \left(1-\exp\left(-\frac{\exp\bigg(-\frac{\beta ^4 k^2}{\left(2 a \beta +\sigma ^2\right)^2}\bigg)}{2 \sqrt{\pi }}\right)\right)dk,
\end{align*}
and 
\begin{align*}
 &\limsup_{N\to\infty}N^{\TailExpon{a}}(\log N)^{\frac{\lambda(a)}{1-\lambda(a)}\frac{\sigma^2}{2\sigma_A^2}}\probability{\MaxQueue{}{}>\TailEvent{a}}\nonumber\\
 \leq& \int_{-\infty}^{\infty}\left(\frac{\sigma_A  \exp \left(-\frac{\beta ^4 k^2 \left(\sigma -\sqrt{2 a \beta +\sigma ^2}\right)^2}{\sigma_A ^2 \left(2 a \beta +\sigma ^2\right)^2}\right)}{2 \sqrt{\pi } \left(\sqrt{2 a \beta +\sigma ^2}-\sigma
   \right)}+\frac{\sigma_A  \left(\sigma ^2+\sigma_A ^2\right) \exp \left(-\frac{2 \beta ^4 k^2 \left(\sigma ^2 \left(\sqrt{2 a \beta +\sigma ^2}-\sigma \right)+a \beta  \left(\sqrt{2 a \beta +\sigma ^2}-2 \sigma
   \right)\right)}{\sigma_A ^2 \left(2 a \beta +\sigma ^2\right)^{5/2}}\right)}{2 \sqrt{\pi } \sigma  \left(\sigma  \left(\sigma -\sqrt{2 a \beta +\sigma ^2}\right)+\sigma_A ^2\right)}\right)\nonumber\\
    &\cdot \frac{\beta ^2 e^{-\frac{\beta ^4 k^2}{\left(2 a \beta +\sigma ^2\right)^2}}}{\sqrt{\pi } \left(2 a \beta +\sigma
   ^2\right)}dk+1.   
\end{align*}

\end{remark}
\subsection*{Acknowledgments}
This work is part of the research program Complexity in high-tech manufacturing, (partly) financed by the Dutch Research Council (NWO) through contract 438.16.121. 

\newpage
\pagestyle{plain}
\bibliographystyle{plain}

\bibliography{refs2}
\end{document}